\documentclass[11pt,a4paper,leqno,twoside]{amsart}
\usepackage{latexsym,amssymb,amsmath}
\input xy
\xyoption{all}

\def\CC{\mathbb C}

\def\RR{\mathbb R}
\def\HH{\mathbb H}
\def\AA{{\mathbb A}}

\def\OO{\mathbb O}

\def\11{\mathbf 1}
\def\PP{\mathbb P}

\def\e1{\varepsilon_1}
\def\e2{\varepsilon_2}
\def\e3{\varepsilon_3}

\def\P2{{\PP}^2}

\def\00{\underline{0}}
\def\J0{{\cal J}_3(\underline{0})}

\def\PJ0{\PP({\cal J}_3(\underline{0}))}

\def\e{\varepsilon}

\def\AP2{{\AA\PP}^2}
\def\RP2{{\RR\PP}^2}
\def\CP2{{\CC\PP}^2}
\def\HP2{{\HH\PP}^2}
\def\OP2{{\OO\PP}^2}

\newtheorem{theo}{Theorem}[section]
\newtheorem{coro}[theo]{Corollary}
\newtheorem{lemm}[theo]{Lemma}
\newtheorem{prop}[theo]{Proposition}
\newtheorem{conj}[theo]{Conjecture}

\newtheorem{prob}[theo]{Problem}
\newtheorem{ques}[theo] {Questions}

\theoremstyle{definition}

\theoremstyle{remark}
\newtheorem{rema}[theo]{Remark}

\begin{document}
\title[Algebraic central division algebras]{On algebraic central division algebras over
Henselian fields of finite absolute Brauer $p$-dimensions and
residually arithmetic type} 

\author{I.D. Chipchakov}
\address{Institute of Mathematics and Informatics\\Bulgarian Academy
of Sciences\\1113 Sofia, Bulgaria: E-mail address:
chipchak@math.bas.bg}

\begin{abstract}
Let $(K, v)$ be a Henselian field with a residue field $\widehat
K$ and a value group $v(K)$, and let $\mathbb{P}$ be the set of
prime numbers. This paper finds conditions on $K$, $v(K)$ and
$\widehat K$ under which every algebraic associative central
division $K$-algebra $R$ contains a central $K$-subalgebra
$\widetilde R$ decomposable into a tensor product of central
$K$-subalgebras $R _{p}$, $p \in \mathbb{P}$, of finite
$p$-primary dimensions $[R _{p}\colon K]$, such that each
finite-dimensional $K$-subalgebra $\Delta $ of $R$ is isomorphic
to a $K$-subalgebra $\widetilde \Delta $ of $\widetilde R $.
\end{abstract}

\maketitle

\section{\bf Introduction}

\medskip
All algebras considered in this paper are assumed to be
associative with a unit (except in the setting of Proposition
\ref{prop2.5} below). Let $E$ be a field, $E _{\rm sep}$ its
separable closure, Fe$(E)$ the set of finite extensions of $E$ in
$E _{\rm sep}$, $\mathbb{P}$ the set of prime numbers,
$\mathbb{P}_{E} = \{p \in \mathbb{P}\colon p \neq {\rm
char}(E)\}$, and for each $p \in \mathbb{P}$, let $E(p)$ be the
maximal $p$-extension of $E$ in $E _{\rm sep}$, i.e. the
compositum of all finite Galois extensions of $E$ in $E _{\rm
sep}$ whose Galois groups are $p$-groups. By the Wedderburn-Artin
theorem (cf. \cite[Theorem~2.1.6]{He}), an Artinian $E$-algebra
$\mathcal{A}$ is simple if and only if it is isomorphic to the
ring $M _{n}(\mathcal{D})$ of $n \times n$-matrices with entries
in a division $E$-algebra $\mathcal{D}$. When this holds, $n$ is
uniquely determined by $\mathcal{A}$, and so is $\mathcal{D}$,
up-to isomorphism; $\mathcal{D}$ is called an underlying division
$E$-algebra of $\mathcal{A}$. The $E$-algebras $\mathcal{A}$ and
$\mathcal{D}$ share a common centre $Z(\mathcal{A})$; we say that
$\mathcal{A}$ is a central $E$-algebra if $Z(\mathcal{A}) = E$.
\par
Denote by Br$(E)$ the Brauer group of $E$, by $s(E)$ the class of
finite-dimensional central simple algebras over $E$, and by $d(E)$
the subclass of division algebras $D \in s(E)$. For each $A \in
s(E)$, let $[A]$ be the equivalence class of $A$ in Br$(E)$, and
let deg$(A)$, ind$(A)$ and exp$(A)$ be the degree, the Schur index
and the exponent of $A$ (the order of $[A]$ as an element of
Br$(E)$), respectively. It is well-known (cf. \cite[Sect.
14.4]{P}) that exp$(A)$ divides ind$(A)$ and shares with it the
same set of prime divisors; also, ind$(A) \mid {\rm deg}(A)$, and
deg$(A) = {\rm ind}(A)$ if and only if $A \in d(E)$. Note that
ind$(B _{1} \otimes _{E} B _{2}) = {\rm ind}(B _{1}){\rm ind}(B
_{2})$ whenever $B _{1}, B _{2} \in s(E)$ and g.c.d.$\{{\rm ind}(B
_{1}), {\rm ind}(B _{2})\} = 1$; equivalently, $B _{1} ^{\prime }
\otimes _{E} B _{2} ^{\prime } \in d(E)$, if $B _{j} ^{\prime }
\in d(E)$, $j = 1, 2$, and g.c.d.$\{{\rm deg}(B _{1} ^{\prime }),
{\rm deg}(B _{2} ^{\prime })\}$ $= 1$ (see \cite[Sect. 13.4]{P}).
Since Br$(E)$ is an abelian torsion group, and ind$(A)$, exp$(A)$
are invariants both of $A$ and the class $[A]$, these results
indicate that Br$(E)$ decomposes into the direct sum of its
$p$-components Br$(E) _{p}\colon p \in \mathbb{P}$, and the study
of the restrictions on the pairs ind$(A)$, exp$(A)$, $A \in s(E)$,
reduces to the special case of $p$-power pairs, for an arbitrary
fixed $p \in \mathbb{P}$. Also, they are crucial for the proof of
Brauer's primary tensor product decomposition theorem for algebras
$D \in d(E)$ (cf. \cite[Sect. 14.4]{P}).
\par
\medskip
This paper is essentially a continuation of \cite{Ch2}. It shows
that if $K _{m}$ is an $m$-local field with a finite residue field
$K _{0}$, for some integer $m \ge 2$, then every algebraic central
division $K _{m}$-algebra $R$ possesses a subalgebra $\widetilde
R$ decomposable into the tensor product $\otimes _{p \in
\mathbb{P}} R _{p}$ of $K _{m}$-algebras $R _{p} \in d(K _{m})$ of
$p$-power degrees, which contains as a $K _{m}$-subalgebra an
isomorphic copy of each finite-dimensional $K _{m}$-subalgebra of
$R$. When $m = 1$, the stated result is contained in
\cite[Theorem~4.1]{Ch2}. It raises interest in two well-known open
problems. The first one, posed by Kurosh in \cite{Ku}, is the
problem of finding whether an algebraic central division algebra
$Y$ over a field $E$ is necessarily locally finite-dimensional
(abbr., LFD), that is, whether every finite subset $\mathcal{S}$
of $Y$ generates a finite-dimensional $E$-subalgebra $Y
_{\mathcal{S}}$ of $Y$. The second one, posed by Koenigsmann in
\cite{Koe}, asks whether the maximal prosolvable extension $F
_{\rm sol}$ of a perfect field $F$ in $F _{\rm sep}$ is a field of
dimension dim$(F _{\rm sol}) \le 1$, i.e. Br$(F _{\rm sol}
^{\prime }) = \{0\}$, for every finite extension $F _{\rm sol}
^{\prime }/F _{\rm sol}$ (see Remark \ref{rema5.2}). As usual, a
Galois extension $E ^{\prime }/E$ is called prosolvable if the
Galois group $\mathcal{G}(E ^{\prime }/E)$ is prosolvable, i.e.
$\mathcal{G}(E ^{\prime }/E)$ is a profinite groups whose finite
homomorphic images are solvable groups.
\par
By a $1$-local field, we mean a complete discrete valued field,
and when $m \ge 2$, an $m$-local field with an $m$-th residue
field $K _{0}$ means a complete field $K _{m}$ with respect to a
discrete valuation $v _{m}$, such that the residue field $K
_{m-1}$ of $(K _{m}, v _{m})$ is an $(m - 1)$-local field with an
$(m - 1)$-th residue field $K _{0}$. If $m \ge 2$ and $\omega
_{m-1}$ is the standard $\mathbb{Z} ^{m-1}$-valued valuation of $K
_{m-1}$, then the composite valuation $\omega _{m} = \omega _{m-1}
{\ast }v _{m}$ is the standard $\mathbb{Z} ^{m}$-valued valuation
of $K _{m}$. It is known that $\omega _{m}$ is Henselian (cf.
\cite[Proposition~A.15]{TW}) and $K _{0}$ equals the residue field
of $(K _{m}, \omega _{m})$. Here $\mathbb{Z} ^{m}$ is viewed as an
ordered abelian group by the inverse-lexicographic ordering. Our
main result on algebraic central division $K _{m}$-algebras
concerns the case where $K _{0}$ is a virtually perfect field of
arithmetic type (see Definition~1 and page \ref{pagerefdefi1}).
\par
\medskip
\section{Background and motivation}
\par
\medskip
The study of the influence of a field $E$ upon the sets of
index-exponent $p$-power pairs ind$(A _{p}), {\rm exp}(A
_{p})\colon A _{p} \in s(E), [A _{p}] \in {\rm Br}(E) _{p}, p \in
\mathbb{P}$, leads naturally to the notion of Brauer
$p$-dimensions Brd$_{p}(E)\colon p \in \mathbb{P}$, introduced in
\cite{ABGV}, for each $p$, as follows: we say that Brd$_{p}(E) = n
_{p} < \infty $ if $n _{p}$ is the least integer $\ge 0$, for
which ind$(A _{p}) \mid {\rm exp}(A _{p}) ^{n_{p}}$ whenever $A
_{p} \in s(E)$ and $[A _{p}] \in {\rm Br}(E) _{p}$; if no such $n
_{p}$ exists, we put Brd$_{p}(E) = \infty $. For instance,
Brd$_{p}(E) \le 1$, for all $p \in \mathbb P$, if and only if
deg$(D) = {\rm exp}(D)$, for each $D \in d(E)$; Brd$_{p'}(E) = 0$,
for some $p' \in \mathbb P$, if and only if Br$(E) _{p'}$ is
trivial.
\par
The absolute Brauer $p$-dimension of $E$ is defined to be the
supremum abrd$_{p}(E)$ of Brd$_{p}(R)\colon R \in {\rm Fe}(E)$. It
is known that $m$-local fields $E _{m}$ with finite $m$-th residue
fields $E _{0}$ have finite absolute Brauer $p$-dimensions, for
each $p \in \mathbb{P}$ (this result is mainly due to Khalin, see
\cite[Propositions~3.4, 3.5 and pp.~318-319]{Ch6}, for more
details and further references); formulae for abrd$_{p}(E _{m})$,
$p \in \mathbb{P}_{E _{0}}$, are also contained in
\cite[Sect.~7]{Br} and \cite[Proposition~5.4 and (3.3)]{Ch7}).
Note further that, by Albert-Hochschild's theorem (cf. \cite[Ch.
II, 2.2]{S1}), abrd$_{p}(E) = 0$, $p \in \mathbb{P}$, if and only
if $E$ is a field of dimension dim$(E) \le 1$, i.e. Br$(R) =
\{0\}$, for every finite extension $R/E$. When $E$ is a perfect
field, we have dim$(E) \le 1$ if and only if the absolute Galois
group $\mathcal{G}_{E} = \mathcal{G}(E _{\rm sep}/E)$ is a
projective profinite group, in the sense of \cite{S1}; it is known
that a profinite group $G$ is projective if and only if its
cohomological dimension cd$(G)$ is $\le 1$ (equivalently, the
Sylow pro-$p$ subgroups $G _{p}$ of $G$ are free pro-$p$ groups
whenever $p \in \mathbb{P}$ and $G _{p} \neq \{1\}$, see \cite[Ch.
I, 5.9, Corollary~2]{S1}). Also, class field theory shows that
Brd$_{p}(E) = {\rm abrd}_{p}(E) = 1$ and dim$(E(p)) \le 1$, for
all $p \in \mathbb{P}$, if $E$ is a global or $1$-local field with
a finite residue field; hence, $E$ is a field of arithmetic type,
in the following sense:
\par
\medskip
{\bf Definition 1.} A field $K$ is said to be of arithmetic type,
if abrd$_{p}(K)$ is finite and abrd$_{p}(K(p)) = 0$, for each $p
\in \mathbb{P}$.
\par
\medskip
\label{pagerefdefi1} It is easily verified that if $E$ is a field
of characteristic $q$, then the inequality Brd$_{p}(E ^{\prime })
\le {\rm abrd}_{p}(E)$ holds whenever $p \in \mathbb{P}_{E}$ and
$E ^{\prime }/E$ is an algebraic extension (see
\cite[Lemma~3.7~(a)]{Ch9} and \cite[(1.2)]{Ch3}). When $q > 0$, we
consider almost exclusively the special case where $E$ is a
virtually perfect field, i.e. $E$ is a finite extension of its
subfield $E ^{q} = \{\lambda ^{q}\colon \lambda \in E\}$. The
class of virtually perfect fields is closed under taking algebraic
extensions; this is implied by the former part of the following
facts concerning finite extensions $E ^{\prime }/E$:
\par
\medskip\noindent
(2.1) (a) $[E ^{\prime }\colon E ^{\prime q}] = [E\colon E ^{q}]$
(cf. \cite[Lemma~2.12]{BH} and \cite[Ch. V, \S{6}]{L});
\par
(b) Brd$_{q}(E ^{\prime }) \le \kappa $, provided that $[E\colon E
^{q}] = q ^{\kappa }$, where $\kappa \in \mathbb{N}$ (for a proof,
see \cite[Lemma~3.7]{Ch9}).
\par
\medskip\noindent
In order to present our research in a unified manner, we say that
$E$ is virtually perfect also in the case of $q = 0$. The class of
virtually perfect fields of finite absolute Brauer $p$-dimensions,
for all $p \in \mathbb{P}$, is singled out by the strong
restrictions on the structure of their central division
LFD-algebras, imposed by \cite[Theorem~3.1]{Ch9}, as follows:
\par
\smallskip
\begin{prop}
\label{prop2.1} Let $R$ be a central division {\rm LFD}-algebra
over a virtually perfect field $K$ with {\rm abrd}$_{p}(K)$
finite, for all $p \in \mathbb{P}$. Then there exist integers
$k(p) \ge 0\colon p \in \mathbb{P}$, and a central $K$-subalgebra
$\widetilde R$ of $R$, such that:
\par
{\rm (a)} $\widetilde R$ is $K$-isomorphic to the tensor product
$\otimes _{p \in \mathbb P} R _{p}$, where $\otimes = \otimes
_{K}$ and $R _{p} \in d(K)$ is a $K$-subalgebra of $R$ of degree
$p ^{k(p)}$, for each $p$;
\par
{\rm (b)} Every finite-dimensional $K$-subalgebra $\Delta $ of $R$
is isomorphic to a subalgebra of $\widetilde R$; hence, $R$
possesses a subalgebra $\Sigma _{\Delta } \in d(K)$ which includes
a $K$-isomorphic copy of $\Delta $;
\par
{\rm (c)} For each $p \in \mathbb{P}$, $k(p)$ is the greatest
integer for which there is $r _{p} \in R$ with $[K(r _{p})\colon
K]$ divisible by $p ^{k(p)}$; $r _{p}$ can be chosen to be
separable over $K$;
\par\noindent
Furthermore, if $[R\colon K]$ is countable, then $\widetilde R$ is
isomorphic to $R$.
\end{prop}
\par
\medskip
The proof of Proposition \ref{prop2.1}, presented in
\cite[Sect.~3]{Ch9}, has been obtained by applying basic results on
separable field extensions together with classical theorems on
Artinian central simple algebras (see \cite{L} and \cite{He}). The
following lemma (proved as \cite[Lemma~8.2]{Ch2}) plays a crucial
role in this proof.
\par
\smallskip
\begin{lemm}
\label{lemm2.2} Let $D$ be a finite-dimensional simple algebra
over a field $K$. Suppose that the centre $B$ of $D$ is a
compositum of extensions $B _{1}$ and $B _{2}$ of $K$ of
relatively prime degrees, and the following conditions are
fulfilled:
\par
{\rm (a)} $[D\colon B] = n ^{2}$ and $D$ possesses a maximal
subfield $E$ such that \par\noindent $[E\colon B] = n$ and $E =
B\widetilde E$, for some separable extension $\widetilde E/K$ of
degree $n$;
\par
{\rm (b)} $p > n$, for every $p \in \mathbb P$ dividing $[B\colon
K]$;
\par
{\rm (c)} $D \cong D _{i} \otimes _{B _{i}} B$ as a $B$-algebra, for
some $D _{i} \in s(B _{i})$, $i = 1, 2$.
\par\noindent
Then there exist $\widetilde D \in s(K)$ with $[\widetilde D\colon K]
= n ^{2}$, and isomorphisms of
\par\noindent
$B _{i}$-algebras $\widetilde D \otimes _{K} B _{i} \cong D _{i}$, $i
= 1, 2, 3$, where $B _{3} = B$ and $D _{3} = D$.
\end{lemm}
\par
\medskip
By Theorem~3.1 of \cite{Ch9}, the conclusion of Proposition
\ref{prop2.1} remains valid if $K$ is a field with char$(K) = q$,
abrd$_{p}(K) < \infty $, $p \in \mathbb{P}_{K}$, and in case $q >
0$, there exists $\mu \in \mathbb{N}$, such that Brd$_{q}(K
^{\prime }) \le \mu $, for every finite extension $K ^{\prime
}/K$. By (2.1) (b), the latter condition is fulfilled if $q > 0$
and $[K\colon K ^{q}] = q ^{\mu }$; also, the existence of $\mu $
is sometimes possible when $[K\colon K ^{q}] = \infty $ and
abrd$_{p}(K)$ are finite, for all $p \in \mathbb{P}_{K}$. More
precisely, for each pair $q \in \mathbb{P}$, $n \in \mathbb{N}$,
there is a field $E _{n}$ with the following properties (see
Proposition \ref{prop4.9}):
\par
\medskip\noindent
(2.2) char$(E _{n}) = q$, $[E _{n}\colon E _{n} ^{q}] = \infty $,
Brd$_{p}(E _{n}) = {\rm abrd}_{p}(E _{n}) = [n/2]$, for all $p \in
\mathbb{P}_{E_{n}}$, and Brd$_{q}(E _{n} ^{\prime }) = n - 1$, for
every finite extension $E _{n} ^{\prime }/E _{n}$.
\par
\medskip\noindent
It should be pointed out that if $E$ is a field with char$(E) = q
> 0$ and $F/E$ is a finitely-generated extension of transcendence
degree $\nu > 0$, then Brd$_{q}(F) < \infty $ if and only if
$[E\colon E ^{q}] < \infty $ \cite[Theorem~2.2]{Ch5}; when
$[E\colon E ^{q}] = q ^{u}$, for some integer $u \ge 0$, we have
$[F\colon F ^{q}] = q ^{u+\nu }$, which means that abrd$_{q}(F)
\le \nu + u$. This attracts interest in the following open
problem:
\par
\medskip
\begin{prob}
\label{prob2.2} Assume that $F/E$ is a finitely-generated field
extension, where {\rm abrd}$_{p}(E) < \infty $, for some $p \in
\mathbb{P}_{E}$. Find whether {\rm abrd}$_{p}(F) < \infty $.
\end{prob}
\par
\medskip
Global fields satisfy the conditions of Proposition \ref{prop2.1},
since they are virtually perfect of absolute Brauer $p$-dimensions
$1$, for all $p \in \mathbb{P}$ (cf. \cite[Example~4.1.3]{Ef}). As
shown more recently by Matzri \cite{Mat}, Proposition
\ref{prop2.1} also applies to any field $K$ of finite Diophantine
dimension, that is, to any field $K$ of type $C _{m}$, for some
integer $m \ge 0$. By type $C _{m}$ (or a $C _{m}$-field, in the
sense of Lang), we mean that every nonzero homogeneous polynomial
$f[X _{1}, \dots , n] \in K[X _{1}, \dots , X _{n}]$ of degree
$d$, where $n > d ^{m}$, has a nontrivial zero over $K$. For
example, algebraically closed fields are of type $C _{0}$; finite
fields are of type $C _{1}$, by the Chevalley-Warning theorem (cf.
\cite[Theorem~6.2.6]{GiSz}), and by Lang's theorem, so are
complete discrete valued fields with algebraically closed residue
fields (see \cite[Ch. II, 3.3]{S1}). The class of $C _{1}$-fields
contains every pseudo algebraically closed (abbr., PAC) field of
characteristic zero (cf. \cite{Kol} and \cite[Remark~21.3.7]{FJ});
in characteristic $q > 0$, perfect PAC fields are of type $C _{2}$
\cite[Theorem~21.3.6]{FJ}. The class of fields of finite
Diophantine dimensions consists of virtually perfect fields and is
closed under the formation of field extensions of finite
transcendence degree (by the Lang-Nagata-Tsen theorem \cite{Na}),
and under taking Laurent formal power series fields in one
variable (see \cite{Gr}). Therefore, the above-noted result of
\cite{Mat} significantly expands the applicability of
\cite[Theorem~3.1]{Ch9}.
\par
\medskip
Proposition \ref{prop2.1} takes a step towards a theory which is
able to reduce the research into central division LFD-algebras
over finitely-generated extensions $E$ of fields $E _{0}$ with
certain specific properties to the study of $E$-algebras $D \in
d(E)$ (some results of this kind can be found, e.g., in
\cite[Sect. 5]{Ch2}). An affirmative solution to Problem
\ref{prob2.2} would allow to take as $E _{0}$ any virtually
perfect field of finite absolute Brauer $p$-dimensions, for all $p
\in \mathbb{P}$. At the same time, Proposition \ref{prop2.1}
attracts interest in the following conjecture:
\par
\medskip
\begin{conj}
\label{conj2.3} Assume that $R$ is an algebraic central division
algebra over a virtually perfect field $K$ with {\rm abrd}$_{p}(K)
< \infty $, for all $p \in \mathbb{P}$. Then $R$ possesses a
$K$-subalgebra $\widetilde R$ subject to the following
restrictions:
\par
{\rm (a)} $\widetilde R$ is $K$-isomorphic to the tensor product
$\otimes _{p \in \mathbb P} R _{p}$, where $\otimes = \otimes
_{K}$ and $R _{p} \in d(K)$ is a $K$-subalgebra of $R$ of
$p$-power degree $p ^{k(p)}$, for each $p \in \mathbb P$;
\par
{\rm (b)} Every $K$-subalgebra $\Delta $ of $R$, which is {\rm
LFD} of at most countable dimension, is embeddable in $\widetilde
R$ as a $K$-subalgebra; also, $K$ equals the centralizer $C
_{R}(\widetilde R) = \{c \in R\colon c\tilde r = \tilde rc, \tilde
r \in \widetilde R\}$, and for each $p \in \mathbb{P}$, $k(p)$ is
the maximal integer for which there is $\rho _{p} \in R$ with
$[K(\rho _{p})\colon K]$ divisible by $p ^{k(p)}$.
\end{conj}
\par
\smallskip
\begin{rema}
\label{rema2.4} It is easily verified that if $K$, $R$ and
$\widetilde R$ are admissible by Conjecture \ref{conj2.3}, and
$[R\colon K]$ is countable but $R$ is not LFD, then the set
$\Sigma _{R}$ of those $K$-subalgebras of $R$, which are
isomorphic to $\widetilde R$, is partially ordered by
set-theoretic inclusion and satisfies the conditions of Zorn's
lemma. Therefore, $\Sigma _{R}$ contains a maximal element
$\widetilde R _{1}$, whence, $\widetilde R _{1} \neq R$ and
$K$-subalgebras $R _{1}$ of $R$ properly including $\widetilde R
_{1}$ are not LFD. In addition, it is not difficult to see that $R
_{1}$ has infinite dimension over its centre $Z(R _{1})$ (cf.
\cite[(1.3)]{Ch3}). When $R$ is a finitely-generated $K$-algebra,
it has a set of generators $\omega _{R} = \{\omega _{0}, \omega
_{1}, \dots , \omega _{m^{2}}\}$, such that $\omega _{1}, \dots ,
\omega _{m^{2}}$ is a basis of some $K$-subalgebra $Y _{m} \in
d(K)$ of $R$ of degree $m$, for some $m \in \mathbb{N}$.
\end{rema}
\par
\medskip
Besides possible applications to the Kurosh problem for algebraic
central division algebras, one may consider other major aspects of
Conjecture \ref{conj2.3} restricting to the case of $[R\colon K] =
\infty $. For reasons clarified in the sequel, in this paper we
assume further that $R$ belongs to the class of $K$-algebras of
linearly bounded degree, introduced in \cite{Am}, as follows:
\par
\medskip
{\bf Definition 2.} An algebraic algebra $\Psi $ over a field $F$
is said to be an algebra of linearly (or locally) bounded degree
(briefly, an LBD-algebra), if the following condition holds, for
any finite-dimensional $F$-subspace $V$ of $\Psi $: there exists
$n(V) \in \mathbb{N}$, such that $[F(v)\colon F] \le n(V)$, for
each $v \in V$.
\par
\medskip
Clearly, if the Kurosh problem is solved affirmatively, then
Conjecture \ref{conj2.3} will turn out to be a restatement of
Proposition \ref{prop2.1}. The solution will be positive if and
only if so are the answers to the following open questions:
\par
\medskip
\begin{ques}
\label{ques2.4} Let $F$ be a field.
\par
{\rm (a)} Find whether algebraic division $F$-algebras are {\rm
LBD}-algebras over $F$.
\par
{\rm (b)} Find whether division {\rm LBD}-algebras over $F$ are
{\rm LFD}.
\end{ques}
\par
\medskip
It is worth noting that the general Kurosh problem, namely, the
problem of finding whether algebraic algebras over a field are LFD
(also posed in \cite{Ku}) is solved, generally, in the negative.
This has been shown by Golod in \cite{Go} along the way to the
proof of the Golod-Shafarevich theorem (see \cite{GS}). For
convenience of the reader, we prove the following proposition
which shows that algebraic unital algebras over countable fields
need not be LBD, and unital LBD-algebras over non-countable fields
need not be LFD.
\par
\medskip
\begin{prop}
\label{prop2.5} Let $R _{0}$ be an associative nil-algebra over an
infinite field $E$, and $R$ be the $E$-algebra obtained from $R
_{0}$ by adjunction of a unit. Then:
\par
{\rm (a)} $R$ is an algebraic associative $E$-algebra;
\par
{\rm (b)} $R$ is a finitely-generated {\rm LBD}-algebra with
$[R\colon E]$ infinite, in case $E$ is non-countable and $R _{0}$
is a finitely-generated algebra with $[R _{0}\colon E]$ infinite;
\par
{\rm (c)} $R$ is not an {\rm LBD}-algebra over $E$, provided that
the polynomial ring $R _{0}[X]$ in an indeterminate $X$ is not
nil.
\end{prop}
\par
\medskip
The existence of a nil-algebra $R _{0}$ over $E$ satisfying the
conditions of Proposition \ref{prop2.5} (b) has been proved by
Golod in \cite{Go}. The existence of a countable field $E$ and a
nil-algebra $R _{0}$ over $E$, admissible by Proposition
\ref{prop2.5} (c), has been established by Smoktunowicz (see
\cite[Theorem~12]{Sm}).
\par
\medskip
\begin{proof} Proposition \ref{prop2.5} (a) is obtained by
straightforward calculations from the definition of $R$ and the
assumption on $R _{0}$; also, it is easily verified that $[R\colon
E] = \infty $ and $R$ is a finitely-generated $E$-algebra,
provided that so is $R _{0}$ with $[R _{0}\colon E] = \infty $.
This, combined with Proposition \ref{prop2.5} (a) and Amitsur's
theorem (cf. \cite[Theorem~9]{Am}), proves Proposition
\ref{prop2.5} (b). We turn to the proof of Proposition
\ref{prop2.5} (c), so we assume that $R _{0}[X]$ is not a
nil-ring, and show that $R \otimes _{E} E(X)$ is not an algebraic
algebra over the rational function field $E(X)$. Consider a
non-nilpotent element $r(X) = \sum _{i=0} ^{\nu } r _{i}X ^{i}$ of
$R _{0}(X)$, where $r _{i} \in R _{0}$ for each index $i$. It is
not difficult to see that, for each $n \in \mathbb{N}$, there
exists a pair $\alpha _{n} \in E$, $s _{n} \in \mathbb{N}$, such
that $s(n) > n$, $r(\alpha _{n}) ^{s(n)} = 0$ and $r(\alpha _{n})
^{s(n)-1} \neq 0$. This enables one to prove by assuming the
opposite that $r(X)$ is a transcendental element over $E(X)$.
Hence, by Amitsur's lemma (see \cite[Lemma~6]{Am}), $R$ is not an
LBD-algebra over $E$, as claimed.
\end{proof}
\par
\medskip
Being closely related to the Kurosh problem, Questions
\ref{ques2.4} (a) and (b) make interest in their own right. For
example, the main results of \cite{Ch2} show that an affirmative
answer to Question \ref{ques2.4}~(a) would prove Conjecture
\ref{conj2.3} in the special case where $K$ is a virtually perfect
field of arithmetic type. Note also that the answer to Question
\ref{ques2.4}~(b) is positive when $F$ is a noncountable field;
this is a special case of \cite[Theorem~9]{Am}, which asserts that
all algebraic $F$-algebras are LBD (see also
\cite[Theorem~5]{Am}). In addition, it follows from \cite[Lemma~6
and Theorem~9]{Am} that if $A$ is an LBD-algebra over an infinite
field $E$, then so is $A \otimes _{E} E ^{\prime }$ (and, more
generally, the full matrix ring $M _{n}(A \otimes _{E} E ^{\prime
})$, for each $n \in \mathbb{N}$) over any extension $E ^{\prime
}$ of $E$; in view of Proposition \ref{prop2.5} (c) and
\cite[Theorem~12]{Sm}, the conclusion of the stated lemma need not
be true if $A$ is merely algebraic and $E$ is countable. The cited
results of \cite{Am}, combined with Wedderburn's structure theorem
(the finite-dimensional version of the Wedderburn-Artin theorem),
indicate that, for each LBD-algebra $A$ over an infinite field
$E$, any extension $E ^{\prime }$ of $E$, and each $B \in s(E
^{\prime })$, the tensor product $A \otimes _{E} B$ is an
LBD-algebra over $E ^{\prime }$. This fact is repeatedly (and
often, implicitly) used in the present paper for proving
Conjecture \ref{conj2.3}, under the hypothesis that $K$ lies in
some classes of Henselian fields of finite absolute Brauer
$p$-dimensions and with virtually perfect residue fields of
arithmetic type. Taking into account the scope of Amitsur's
theorem, we recall that the class $\mathcal{HNF}$ of Henselian
noncountable fields contains every maximally complete field, i.e.
any nontrivially valued field $(K, v)$ which does not admit a
valued proper extension with the same value group and residue
field. For instance, $\mathcal{HNF}$ contains the generalized
(formal) power series field $K _{0}((\Gamma ))$ over a field $K
_{0}$, where $\Gamma $ is a nontrivial ordered abelian group, and
$v$ is the standard valuation of $K _{0}((\Gamma ))$ trivial on $K
_{0}$ (see \cite[Example~4.2.1 and Theorem~18.4.1]{Ef}). Moreover,
for each $m \in \mathbb{N}$, every $m$-local field is
$\mathcal{HNF}$ with respect to its $\mathbb{Z} ^{m}$-valued
valuation.
\par
\medskip
\section{\bf Statements of the main results}
\par
\medskip
Let $R$ be an algebraic central division algebra over a virtually
perfect field $K$ with abrd$_{p}(K) < \infty $, for all $p \in
\mathbb{P}$. Clearly, if $R$ possesses a $K$-subalgebra
$\widetilde R$ as described by Conjecture \ref{conj2.3}, then
there is a sequence $\bar k = k(p)$, $p \in \mathbb{P}$, of
integers $\ge 0$, such that the numbers $p ^{k(p)+1}\colon p \in
\mathbb{P}$, do not divide $[K(r)\colon K]$, for any $r \in R$.
The existence of $\bar k$ is always guaranteed if $R$ is an
LBD-algebra over $K$ (cf. \cite[Lemma~3.9]{Ch2}). When $k(p) =
k(p) _{R}$ is the minimal integer satisfying the stated condition,
it is called a $p$-power of $R/K$. In this setting, the notion of
a $p$-splitting field of $R/K$ is defined as follows:
\par
\medskip
{\bf Definition 3.} Let $K ^{\prime }$ be a finite extension of $K$,
$R ^{\prime }$ the underlying (central) division $K ^{\prime
}$-algebra of $R \otimes _{K} K ^{\prime }$, and $\gamma (p)$ the
integer singled out by the Wedderburn-Artin $K ^{\prime
}$-isomorphism $R \otimes _{K} K ^{\prime } \cong M _{\gamma (p)}(R
^{\prime })$. We say that $K ^{\prime }$ is a $p$-splitting field of
$R/K$ if $p ^{k(p)}$ divides $\gamma (p)$.
\par
\medskip\noindent
Note that the class of $p$-splitting fields of a central division
LBD-algebra $R$ over a virtually perfect field $K$ with
abrd$_{p'}(K) < \infty $, $p' \in \mathbb{P}$, is closed under the
formation of finite extensions. Indeed, it is known (cf.
\cite[Lemma~4.1.1]{He}) that $R \otimes _{K} K ^{\prime }$ is a
central simple $K ^{\prime }$-algebra, for any field extension $K
^{\prime }/K$. This algebra is a left (and right) vector space
over $R$ of dimension equal to $[K ^{\prime }\colon K]$, which
implies it is Artinian whenever $[K ^{\prime }\colon K]$ is
finite. Since $R \otimes _{K} K _{2}$ and $(R \otimes _{K} K _{1})
\otimes _{K _{1}} K _{2}$ are isomorphic $K _{2}$-algebras, for
any tower of field extensions $K \subseteq K _{1} \subseteq K
_{2}$ (cf. \cite[Sect. 9.4, Corollary~(a)]{P}), these observations
enable one to deduce our assertion about the class of
$p$-splitting fields of $R/K$ from the Wedderburn-Artin theorem
(and basic properties of tensor products of matrix algebras, see
\cite[Sect. 9.3, Corollary~b]{P}). Further results on $k(p)$ and
the $p$-power of the underlying division $K ^{\prime }$-algebra of
$R \otimes _{K} K ^{\prime }$, obtained in case $K ^{\prime }/K$
is a finite extension, are presented in Section~5 (see Lemma
\ref{lemm5.1}). They have been proved in \cite[Sect. 3]{Ch2} under
the extra hypothesis that dim$(K _{\rm sol}) \le 1$, where $K
_{\rm sol}$ is the maximal prosolvable extension of $K$ in $K
_{\rm sep}$. These results partially generalize well-known facts
about algebras $D \in d(E)$ for an arbitrary field $E$, leaving
open the question of whether the condition that dim$(K _{\rm sol})
\le 1$ is superfluous, and whether it is essential for the
validity of the derived information (see Remark \ref{rema5.2}).
\par
\medskip
The purpose of this paper is to prove Conjecture \ref{conj2.3} for
two types of Henselian fields. Our first main result can be stated
as follows:
\par
\medskip
\begin{theo}
\label{theo3.1} Let $K = K _{m}$ be an $m$-local field with a
virtually perfect $m$-th residue field $K _{0}$, for some integer
$m > 0$, and let $R$ be an algebraic central division $K$-algebra.
Assume that {\rm char}$(K _{0}) = q$ and $K _{0}$ is of arithmetic
type. Then $R$ has a $K$-subalgebra $\widetilde R$ with the
properties claimed by Conjecture \ref{conj2.3}.
\end{theo}
\par
\medskip
When char$(K _{m}) = {\rm char}(K _{0})$, $(K _{m}, w _{m})$ is
isomorphic to the $m$-fold iterated formal Laurent power series
field $\mathcal{K} _{m} := K _{0}((X _{1})) \dots ((X _{m}))$
(that is, $K _{m'} = K _{m'-1}((X _{m'}))$, for every index $m' >
0$), considered with its standard $\mathbb{Z} ^{m}$-valued
valuation, say $\tilde w _{m}$, acting trivially on $K _{0}$; in
particular, this holds if char$(K _{0}) = 0$. It is known that
$(\mathcal{K} _{m}, \tilde w _{m})$ is maximally complete (cf.
\cite[Sect. 18.4]{Ef}). As $K _{0}$ is virtually perfect, whence,
so is $\mathcal{K} _{m}$, this enables one to prove the assertion
of Theorem \ref{theo3.1} by applying our next result to
$(\mathcal{K} _{m}, \tilde w _{m})$:
\par
\medskip
\begin{theo}
\label{theo3.2} Let $(K, v)$ be a Henselian field with $\widehat
K$ of arithmetic type, and suppose that {\rm char}$(K) = {\rm
char}(\widehat K) = q$, $K$ is virtually perfect and {\rm
abrd}$_{p}(K)$ is finite, for each $p \in \mathbb{P}_{\widehat
K}$. Then every central division {\rm LBD}-algebra $R$ over $K$
has a central $K$-subalgebra $\widetilde R$ admissible by
Conjecture \ref{conj2.3}.
\end{theo}
\par
\medskip
The assertions of Theorems \ref{theo3.1} and \ref{theo3.2} are
well-known in case $R \in d(K)$. When $[R\colon K] = \infty $,
they can be deduced from Lemma \ref{lemm2.2} and the following
generalization of \cite[Lemma~8.3]{Ch2}, by the method of proving
\cite[Theorem~4.1]{Ch2} (see the end of Section 5).
\par
\medskip
\begin{lemm}
\label{lemm3.3} Assume that $K$ is a field and $R$ is a central
division $K$-algebra, which satisfy the conditions of some of
Theorems \ref{theo3.1} or \ref{theo3.2}. Then, for each $p \in
\mathbb{P}$, $K$ has a finite extension $E _{p}$ in $K(p)$ that is
a $p$-splitting field of $R/K$.
\end{lemm}
\par
\medskip
The conditions of Lemma \ref{lemm3.3} ensure that dim$(K _{\rm
sol}) \le 1$ (see Lemmas \ref{lemm6.4} and \ref{lemm6.5} below);
hence, by Lemma \ref{lemm5.1}, a finite extension $E _{p}$ of $K$
in $K _{\rm sep}$ is a $p$-splitting field of $R/K$, for some $p
\in \mathbb{P}$, if and only if $p \nmid [E _{p}(\rho _{p})\colon
E _{p}]$, for any element $\rho _{p}$ of the underlying division
$E _{p}$-algebra $\mathcal{R}_{p}$ of $R \otimes _{K} E _{p}$.
This plays an essential role in the proofs of Lemma \ref{lemm3.3}
and our main theorems. Note also that, in the setting of
Conjecture \ref{conj2.3}, it is unknown whether there exist
$p$-splitting fields $E _{p}$, $p \in \mathbb{P}$, of $R/K$, such
that $E _{p} \subseteq K(p)$, for each $p$. In view of Proposition
\ref{prop2.1} and \cite[Conjecture~1]{Me} (noted at the end of
\cite[Ch. 15]{P}), and since Questions \ref{ques2.4} are open, the
answer is affirmative in all presently known cases. When $R$ is
LFD and $K$ contains a primitive $p$-th root of unity, for each $p
\in \mathbb{P}_{K}$, such an answer follows from Proposition
\ref{prop2.1}, combined with Albert's theory of $q$-algebras over
fields of characteristic $q > 0$ (cf. \cite[Ch. VII,
Theorem~28]{A1}), and the Merkur'ev-Suslin theorem
\cite[(16.1)]{MS} (see also \cite[Theorem~9.1.4 and Ch. 8]{GiSz},
respectively).
\par
\medskip
The earliest draft of this paper is contained in the manuscript
\cite{Ch1}. Here we extend the scope of results of \cite{Ch1},
relying on the theory of division algebras over Henselian fields,
developed in \cite{JW}, and on the contribution to absolute Brauer
$p$-dimensions made in \cite{Mat}, \cite{PS} and other papers. The
progress of the research in these areas allows us to consider the
topic of the present paper in the desired generality; for example,
Theorem \ref{theo3.1} applies to $m$-local fields with finite
$m$-th residue fields, which are not of arithmetic type, in the
sense of Definition~1, for any $m \ge 2$ (see Lemmas
\ref{lemm4.8}, \ref{lemm7.4} and page \pageref{7.5}).
\par
\smallskip
The basic notation, terminology and conventions kept in this paper
are standard and virtually the same as in \cite{He}, \cite{L},
\cite{P} and \cite{Ch5}. Throughout, Brauer groups, value groups
and ordered abelian groups are written additively, Galois groups
are viewed as profinite with respect to the Krull topology, and by
a profinite group homomorphism, we mean a continuous one. For any
algebra $A$, we consider only subalgebras containing its unit.
Given a field $E$, $E ^{\ast }$ denotes its multiplicative group,
$E ^{\ast n} = \{a ^{n}\colon \ a \in E ^{\ast }\}$, for each $n
\in \mathbb N$, and for any $p \in \mathbb P$, $_{p}{\rm Br}(E)$
stands for the maximal subgroup $\{b _{p} \in {\rm Br}(E)\colon \
pb _{p} = 0\}$ of Br$(E)$ of period dividing $p$. We denote by
$I(E ^{\prime }/E)$ the set of intermediate fields of any field
extension $E ^{\prime }/E$, and by Br$(E ^{\prime }/E)$ the
relative Brauer group of $E ^{\prime }/E$ (the kernel of the
scalar extension map Br$(E) \to {\rm Br}(E ^{\prime })$). In case
char$(E) = q > 0$, we write $[a, b) _{E}$ for the $q$-symbol
$E$-algebra generated by elements $\xi $ and $\eta $, such that
\par\noindent
$\eta \xi = (\xi + 1)\eta $, $\xi ^{q} - \xi = a \in E$ and $\eta
^{q} = b \in E ^{\ast }$.
\par
\smallskip
Here is an overview of the rest of this paper. Section 4 includes
preliminaries on Henselian fields used in the sequel. We also state
two lemmas and a corollary, proved in \cite[Sect.~4]{Ch5}, which
fully characterize Henselian fields of residual characteristic zero
as well as $m$-local fields, and maximally complete
equicharacteristic fields (with an emphasis on generalized formal
power series fields) of finite absolute Brauer $p$-dimensions, for
all $p \in \mathbb{P}$; these fields turn out to be virtually
perfect, whence, admissible by Proposition \ref{prop2.1}. Section 5
contains lemmas providing Galois-theoretic and ring-theoretic
ingredients of the proofs of Lemma \ref{lemm3.3} and our main
results. Most of them are stated here without proofs, since they have
been extracted from \cite{Ch2}. These lemmas allow to borrow ideas
from (and use results of) the theory of Artinian central simple
algebras in our considerations; in the first place, to show that
Theorems \ref{theo3.1} and \ref{theo3.2} can be deduced from Lemma
\ref{lemm3.3}, by the method of proving \cite[Theorem~4.1]{Ch2}.
As to the idea to take Lemma \ref{lemm3.3} as a key point in the
proof of our main results, it comes from Lemmas \ref{lemm4.3} and
\ref{lemm6.3}; by the latter lemma, Henselian fields $(K, v)$ with
char$(K) = q > 0$ satisfy abrd$_{q}(K(q)) = 0$, and so do
Henselian discrete valued (abbr., HDV) fields of residual
characteristic $q$. It is also proved in Section 6 that, in the
setting of Theorem \ref{theo3.1} or \ref{theo3.2}, dim$(K _{\rm
sol}) \le 1$ (see Lemmas \ref{lemm6.4} and \ref{lemm6.5}), which
allows to use in our proof results of Section 5. Section 7
collects valuation-theoretic ingredients of the proof of Lemma
\ref{lemm3.3}, including its tame version, stated as Lemma
\ref{lemm7.6}. Finally, we prove Lemma \ref{lemm3.3}. This is done in
Sections 8 and 9, by adapting to Henselian fields the method of
proving \cite[Lemma~8.3]{Ch2} so as to make it possible to apply
Lemma \ref{lemm2.2}. Also, we rely on Lemmas \ref{lemm7.3},
\ref{lemm7.4}, \ref{lemm7.6}, and on lemmas considered in the
preceding three sections.
\par
\medskip
\section{\bf Preliminaries on Henselian fields and their
finite-dimensional division algebras}
\par
\medskip
Let $K$ be a field with a nontrivial valuation $v$, $O _{v}(K) =
\{a \in K\colon \ v(a) \ge 0\}$ the valuation ring of $(K, v)$,
$\mathcal{M} _{v}(K) = \{\mu \in K\colon \ v(\mu ) > 0\}$ the
maximal ideal of $O _{v}(K)$, $O _{v}(K) ^{\ast } = \{u \in
K\colon \ v(u) = 0\}$ the multiplicative group of $O _{v}(K)$,
$v(K)$ and $\widehat K = O _{v}(K)/\mathcal{M} _{v}(K)$ the value
group and the residue field of $(K, v)$, respectively; put $\nabla
_{0}(K) = \{\alpha \in K\colon \alpha - 1 \in \mathcal{M}
_{v}(K)\}$. The valuation $v$ is called discrete if $v(K)$ is an
infinite cyclic group. We say that $v$ is Henselian if it extends
uniquely, up-to equivalence, to a valuation $v _{L}$ on each
algebraic extension $L$ of $K$. This holds, for example, if $K = K
_{v}$ and $v(K)$ is an ordered subgroup of the additive group
$\mathbb R$ of real numbers (cf. \cite[Ch. XII]{L}). Maximally
complete fields are also Henselian, since Henselizations of valued
fields are their immediate extensions (see, e.g.,
\cite[Proposition~15.3.7]{Ef}, or \cite[Corollary~A.28]{TW}). In
order that $v$ be Henselian, it is necessary and sufficient that
any of the following two equivalent conditions is fulfilled (cf.
\cite[Theorem~18.1.2]{Ef} or \cite[Theorem~A.14]{TW}):
\par
\medskip\noindent
(4.1) (a) Given a polynomial $f(X) \in O _{v}(K) [X]$ and an
element $a \in O _{v}(K)$, such that $2v(f ^{\prime }(a)) <
v(f(a))$, where $f ^{\prime }$ is the formal derivative of $f$,
there is a zero $c \in O _{v}(K)$ of $f$ satisfying the equality
$v(c - a) = v(f(a)/f ^{\prime }(a))$;
\par
(b) For each normal extension $\Omega /K$, $v ^{\prime }(\tau (\mu ))
= v ^{\prime }(\mu )$ whenever  $\mu \in \Omega $, $v ^{\prime }$ is
a valuation of $\Omega $ extending $v$, and $\tau $ is a
$K$-automorphism of $\Omega $.
\par
\medskip\noindent
When $v$ is Henselian, so is $v _{L}$, for any algebraic field
extension $L/K$. In this case, we put $O _{v}(L) = O _{v
_{L}}(L)$, $\mathcal{M} _{v}(L) = \mathcal{M} _{v_{L}}(L)$, $v(L)
= v _{L}(L)$, and denote by $\widehat L$ the residue field of $(L,
v _{L})$. Clearly, $\widehat L/\widehat K$ is an algebraic
extension and $v(K)$ is an ordered subgroup of $v(L)$; the index
$e(L/K)$ of $v(K)$ in $v(L)$ is called a ramification index of
$L/K$. By Ostrowski's theorem (see \cite[Sects. 17.1 and
17.2]{Ef}) if $[L\colon K]$ is finite, then $[\widehat L\colon
\widehat K]e(L/K)$ divides $[L\colon K]$, and the integer
$[L\colon K][\widehat L\colon \widehat K] ^{-1}e(L/K) ^{-1}$ is
not divisible by any $p \in \mathbb{P}_{\widehat K}$; also, $L/K$
is defectless, i.e. $[L\colon K] = [\widehat L\colon \widehat
K]e(L/K)$, in the following three cases:
\par
\medskip\noindent
(4.2) (a) If char$(\widehat K) \nmid [L\colon K]$ (apply
Ostrowski's theorem);
\par
(b) If $(K, v)$ is HDV and $L/K$ is separable (see
\cite[Theorem~A.12]{TW});
\par
(c) When $(K, v)$ is maximally complete (cf.
\cite[Theorem~31.22]{Wa}).
\par
\medskip\noindent
Assume that $(K, v)$ is a Henselian field and $\mathcal{I}/K$ is a
finite extension. We say that $\mathcal{I}/K$ is inertial, if
$[\mathcal{I}\colon K] = [\widehat {\mathcal{I}}\colon \widehat
K]$ and $\widehat {\mathcal{I}}/\widehat K$ is a separable
extension; $\mathcal{I}/K$ is said to be totally ramified, if
$e(\mathcal{I}/K) = [\mathcal{I}\colon K]$. Inertial extensions of
$K$ have the following useful properties (see \cite{TW},
Theorem~A.23, Proposition~A.17, and Corollary~A.25):
\par
\medskip
\begin{lemm}
\label{lemm4.1} Let $(K, v)$ be a Henselian field. Then:
\par
{\rm (a)} An inertial extension $\mathcal{I}/K$ is Galois if and
only if $\widehat {\mathcal{I}}/\widehat K$ is Galois. When this
holds, the Galois groups $\mathcal{G}(\mathcal{I}/K)$ and
$\mathcal{G}(\widehat {\mathcal{I}}/\widehat K)$ are isomorphic.
\par
{\rm (b)} The compositum $K _{\rm ur}$ of inertial extensions of
$K$ in $K _{\rm sep}$ is a Galois extension of $K$ with
$\mathcal{G}(K _{\rm ur}/K) \cong \mathcal{G}_{\widehat K}$.
\par
{\rm (c)} Finite extensions of $K$ in $K _{\rm ur}$ are inertial,
and the natural mapping of $I(K _{\rm ur}/K)$ into $I(\widehat K
_{\rm sep}/\widehat K)$ is bijective.
\par
{\rm (d)} For each $K _{1} \in {\rm Fe}(K)$, the intersection $K _{0}
= K _{1} \cap K _{\rm ur}$ equals the maximal inertial extension of
$K$ in $K _{1}$; in addition, $\widehat K _{0} = \widehat K _{1}$.
\end{lemm}
\par
\medskip\noindent
When $(K, v)$ is Henselian, a finite extension $L/K$ is called
tamely ramified if char$(\widehat K) \nmid e(L/K)$, $\widehat
L/\widehat K$ is separable, and $L/K$ is defectless. The next
lemma gives an account of some basic properties of tamely ramified
extensions of $K$ in $K _{\rm sep}$ (see (4.1) (b) and
\cite[Theorems~A.9 (i),(ii) and A.24]{TW}):
\par
\medskip
\begin{lemm}
\label{lemm4.2} Let $(K, v)$ be a Henselian field with {\rm
char}$(\widehat K) = q$, $K _{\rm tr}$ the compositum of tamely
ramified extensions of $K$ in $K _{\rm sep}$, and for each $p \in
\mathbb{P}_{\widehat K}$, let $\varepsilon _{p}$, $\hat
\varepsilon _{p}$ be primitive $p$-th root of unity in $K _{\rm
sep}$ and $\widehat K _{\rm sep}$, respectively. Then:
\par\vskip0.04truecm
{\rm (a)} $K _{\rm tr}/K$ is a Galois extension with
$\mathcal{G}(K _{\rm tr}/K _{\rm ur})$ abelian, and all finite
extensions of $K$ in $K _{\rm tr}$ are tamely ramified.
\par
{\rm (b)} There is $T(K) \in I(K _{\rm tr}/K)$ with $T(K) \cap K
_{\rm ur} = K$ and $T(K).K _{\rm ur} = K _{\rm tr}$; hence, finite
extensions of $K$ in $T(K)$ are tamely and totally ramified.
\par
{\rm (c)} The field $T(K)$ singled out in {\rm (b)} is isomorphic
as a $K$-algebra to $\otimes _{p \in \mathbb{P}_{\widehat K}} T
_{p}(K)$, where $\otimes = \otimes _{K}$, and for each $p$, $T
_{p}(K) \in I(T(K)/K)$ and every finite extension of $K$ in $T
_{p}(K)$ is of $p$-power degree; in particular, $T(K)$ equals the
compositum of the fields $T _{p}(K)$, $p \in \mathbb{P}_{\widehat
K}$.
\par
{\rm (d)} With notation being as in {\rm (c)}, $T _{p}(K) \neq K$,
for some $p \in \mathbb{P}_{\widehat K}$, if and only if $v(K)
\neq pv(K)$; when this holds, $T _{p}(K) \in I(K(p)/K)$ if and
only if $\hat \varepsilon _{p} \in \widehat K$ (equivalently, if
and only if $\varepsilon _{p} \in K$).
\end{lemm}
\par
\medskip
The Henselian property of $(K, v)$ guarantees that $v$ extends to
a unique, up-to equivalence, valuation $v _{D}$ on each $D \in
d(K)$ (cf. \cite[Sect. 1.2.2]{TW}). Put $v(D) = v _{D}(D)$ and
denote by $\widehat D$ the residue division ring of $(D, v _{D})$.
It is known that $\widehat D$ is a division $\widehat K$-algebra,
$v(D)$ is an ordered abelian group and $v(K)$ is an ordered
subgroup of $v(D)$ of finite index $e(D/K)$ (called the
ramification index of $D/K$). Note also that $[\widehat D\colon
\widehat K] < \infty $, and by Ostrowski-Draxl's theorem (cf.
\cite{Dr2} and \cite[Propositions~4.20 and 4.21]{TW}), $[\widehat
D\colon \widehat K]e(D/K) \mid [D\colon K]$ and $[D\colon
K][\widehat D\colon \widehat K] ^{-1}e(D/K) ^{-1}$ has no divisor
\par\vskip0.031truecm\noindent
$p \in \mathbb{P}_{\widehat K}$. The division $K$-algebra $D$ is
said to be inertial if $[D\colon K] = [\widehat D\colon \widehat
K]$; it is called totally ramified if $[D\colon K] = e(D/K)$. We
say that $D/K$ is defectless, if $[D\colon K] = [\widehat D\colon
\widehat K]e(D/K)$; this holds in the following two cases:
\par\medskip\noindent
(4.3) (a) If char$(\widehat K) \nmid [D\colon K]$ (apply the
Ostrowski-Draxl theorem);
\par
(b) If $(K, v)$ is an HDV-field (see \cite[Proposition~2.2]{TY}).
\par
\medskip\noindent
The $K$-algebra $D$ is called nicely semi-ramified (abbr., NSR),
in the sense of \cite{JW}, if $e(D/K) = [\widehat D\colon \widehat
K] = {\rm deg}(D)$, $\widehat D/\widehat K$ is a separable field
extension, and $D$ possesses a maximal subfield which is a totally
ramified extension of $K$ (see \cite[page~148]{JW} and
\cite[Theorem~6.4]{Moun}). It is known that if $D$ is NSR, then
$\widehat D/\widehat K$ is a Galois extension,
$\mathcal{G}(\widehat D/\widehat K)$ is isomorphic to the quotient
group $v(D)/v(K)$, and $D$ decomposes into a tensor product of
cyclic NSR-algebras over $K$ (see \cite[Example~4.3 and
Theorem~4.4]{JW} or \cite[Propositions~8.40, 8.41]{TW}). This
result allows to prove the following lemma:
\par
\medskip
\begin{lemm}
\label{lemm4.3} Let $(K, v)$ be a Henselian field with {\rm
abrd}$_{p}(\widehat K(p)) = 0$, for some $p \in
\mathbb{P}_{\widehat K}$. Then every $\Delta _{p} \in d(K)$ of
$p$-power degree has a splitting field that is a finite extension
of $K$ in $K(p)$.
\end{lemm}
\par
\medskip
Lemma \ref{lemm4.3} shows that if $R$ is a central division
LBD-algebra over a field $K$ satisfying the conditions of some of
Theorems \ref{theo3.1} and \ref{theo3.2}, and if there is a
$K$-subalgebra $\widetilde R$ of $R$ with the properties claimed
by Conjecture \ref{conj2.3}, then for each $p \in \mathbb{P}$ with
at most one exception, $K$ has a finite extension $E _{p}$ in
$K(p)$ that is a $p$-splitting field of $R/K$. This leads to the
idea of using Lemma \ref{lemm3.3} as a basis for the proofs of
Theorems \ref{theo3.1} and \ref{theo3.2} (for further support of
the idea and a step to its implementation, see Lemma
\ref{lemm6.3}).
\par
\medskip
{\it Proof of Lemma \ref{lemm4.3}.} The assertion is obvious if
$\Delta _{p}$ is an NSR-algebra over $K$, or more generally, if
$\Delta _{p}$ is Brauer equivalent to the tensor product of cyclic
division $K$-algebras of $p$-power degrees. When the $K$-algebra
$\Delta _{p}$ is inertial, we have $\widehat \Delta _{p} \in
d(\widehat K)$ (cf. \cite{JW}, Theorem~2.8), so our conclusion
follows from the fact that abrd$_{p}(\widehat K(p)) = 0$ and
$\widehat {K(p)} = \widehat K(p)$,  which ensures
\par\vskip0.047truecm\noindent
that $[\Delta _{p}] \in {\rm Br}(K _{\rm ur} \cap K(p)/K)$. Since,
by \cite[Lemmas~5.14 and 6.2]{JW},
\par\vskip0.056truecm\noindent
$[\Delta _{p}] = [I _{p} \otimes _{K} N _{p} \otimes _{K} T
_{p}]$, for some inertial algebra $I _{p}/K$, an NSR-algebra $N
_{p}/K$, and a tensor product $T _{p}$ of totally ramified cyclic
division $K$-algebras,
\par\vskip0.047truecm\noindent
such that $[I _{p}], [N _{p}]$ and $[T _{p}] \in {\rm Br}(K)
_{p}$, these observations prove Lemma \ref{lemm4.3}.
\par
\medskip
The following lemma characterizes all Henselian fields $(K, v)$
with \par\vskip0.04truecm\noindent abrd$_{p}(K) < \infty $, for
some $p \in \mathbb{P}_{\widehat K}$. Since abrd$_{p}(\widehat K)
\le {\rm abrd}_{p}(K)$ (by \cite[Theorem~2.8]{JW} and
\cite[Theorem~A.23]{TW}), it can be deduced from
\cite[Proposition~6.1, Theorem~5.9 and Remark~6.2]{Ch7} (or
\cite[(3.3) and Theorem~2.3]{Ch7}).
\par
\medskip
\begin{lemm}
\label{lemm4.4} For a Henselian field $(K, v)$, {\rm
abrd}$_{p}(K)$ is finite, for some $p \in \mathbb{P}_{\widehat
K}$, if and only if so are {\rm abrd}$_{p}(\widehat K)$ and the
quotient group $v(K)/pv(K)$.
\end{lemm}
\par
\medskip
Lemma \ref{lemm4.4} and our next lemma show that a maximally
complete field $(K, v)$ with char$(K) = {\rm char}(\widehat K)$
satisfies abrd$_{p}(K) < \infty $, $p \in \mathbb{P}$, if and only
if $\widehat K$ is virtually perfect and for each $p \in
\mathbb{P}$, abrd$_{p}(\widehat K) < \infty $ and $v(K)/pv(K)$ is
finite. When this holds, $K$ is virtually perfect as well (see
\cite[Lemma~3.2]{Ch7}).
\par
\medskip
\begin{lemm}
\label{lemm4.5} Let $(K, v)$ be a Henselian field with {\rm
char}$(\widehat K) = q > 0$. Then:
\par
{\rm (a)} $[\widehat K\colon \widehat K ^{q}]$ and $v(K)/qv(K)$
are finite, in case {\rm Brd}$_{q}(K) < \infty $;
\par
{\rm (b)} The inequality {\rm abrd}$_{q}(K) < \infty $ holds,
provided that $\widehat K$ is virtually perfect and $v$ is
discrete or the following condition is satisfied:
\par
{\rm (i)} {\rm char}$(K) = q$ and $K$ is virtually perfect; in
particular, this occurs if {\rm char}$(K) = q$, $v(K)/qv(K)$ is
finite and $(K, v)$ is maximally complete.
\end{lemm}
\par
\medskip
For a proof of Lemma \ref{lemm4.5} and of the following corollary,
we refer the reader to \cite[Lemma~4.7]{Ch9} and
\cite[Proposition~4.4]{Ch9}, respectively.
\par
\medskip
\begin{coro}
\label{coro4.6} Let $K _{0}$ be a field and $\Gamma $ a nontrivial
ordered abelian group. Then the generalized power series field $K
= K _{0}((\Gamma ))$ satisfies the inequalities {\rm abrd}$_{p}(K)
< \infty $, $p \in \mathbb{P}$, if and only if $K _{0}$ is
virtually perfect, {\rm abrd}$_{p}(K _{0}) < \infty $, $p \in
\mathbb{P}_{K_{0}}$, and the quotient groups $\Gamma /p\Gamma $
are finite, for all $p \in \mathbb{P}$.
\end{coro}
\par
\medskip
\begin{rema}
\label{rema4.7} Given a field $K _{0}$ and an ordered abelian
group $\Gamma \neq \{0\}$, the standard realizability of the field
$K _{1} = K _{0}((\Gamma ))$ as a maximally complete field (with a
value group $\Gamma $ and a residue field $K _{0}$), used in the
proof of Corollary \ref{coro4.6}, allows to determine the sequence
$(b, a) = {\rm Brd}_{p}(K _{1}), {\rm abrd}_{p}(K _{1})\colon p
\in \mathbb{P}$, in the following two cases: {\rm (i)} $K _{0}$ is
a global or a $1$-local field with a finite residue field (see
\cite[Proposition~5.1]{Ch7} and \cite[Corollary~3.6 and
Sect.~4]{Ch6}, respectively); {\rm (ii)} $K _{0}$ is perfect and
dim$(K _{0}) \le 1$ (see \cite[Proposition~3.5]{Ch6} and
\cite[Propositions~5.3, 5.4]{Ch7}). In both cases, $(b, a)$
depends only on $K _{0}$ and $\Gamma $. Moreover, if $(K, v)$ is
Henselian with $\widehat K = K _{0}$ and $v(K) = \Gamma $, then:
{\rm (a)} $(b, a) = {\rm Brd}_{p}(K), {\rm abrd}_{p}(K)$, $p \in
\mathbb{P}_{K_{0}}$; {\rm (b)} abrd$_{q}(K) = {\rm abrd}_{q}(K
_{1})$ and Brd$_{q}(K) = {\rm Brd}_{q}(K _{1})$, if $(K, v)$ is
maximally complete and char$(K) = q > 0$.
\end{rema}
\par
\medskip
Next we show that abrd$_{p}(K _{m}) < \infty $, $p \in
\mathbb{P}$, if $K _{m}$ is an $m$-local field with $m$-th residue
field admissible by Theorem \ref{theo3.2}.
\par
\medskip
\begin{lemm}
\label{lemm4.8} Let $K _{m}$ be an $m$-local field with $m$-th
residue field $K _{0}$. Then {\rm abrd}$_{p}(K _{m}) < \infty $,
for all $p \in \mathbb{P}$, if and only if $K _{0}$ is virtually
perfect and {\rm abrd}$_{p}(K _{0}) < \infty $, for all $p \in
\mathbb{P}_{K_{0}}$; when this holds, $K _{m}$ is virtually
perfect.
\end{lemm}
\par
\smallskip
\begin{proof}
Clearly, it suffices to consider the case where abrd$_{p}(K
_{m-1}) < \infty $, for all $p \in \mathbb{P}$. Then our assertion
follows from Lemmas \ref{lemm4.4} and \ref{lemm4.5}.
\end{proof}
\par
\medskip
The concluding result of this section proves (2.2) and leads to
the following open question: given a field $E$ with char$(E) = q >
0$, $[E\colon E ^{q}] = \infty $ and abrd$_{q}(E) < \infty $, does
there exist $\mu (E) \in \mathbb{N}$, such that Brd$_{q}(E
^{\prime }) \le \mu (E)$, for all finite extensions $E ^{\prime
}/E$? An affirmative answer would allow to drop in Proposition
\ref{prop2.1} the condition that $K$ is a virtually perfect field.
\par
\medskip
\begin{prop}
\label{prop4.9} Let $F _{0}$  be an algebraically closed field of
characteristic $q > 0$ and $F _{n}$: $n \in \mathbb{N}$, be
extensions of $F _{0}$ defined inductively as follows: when $n >
0$, $F _{n} = F _{n-1}((T _{n}))$ is the formal Laurent power
series field in a variable $T _{n}$ over $F _{n-1}$. Then the
following holds, for each $n \in \mathbb{N}$:
\par
{\rm (a)} $F _{n}$ possesses a subfield $\Lambda _{n}$ that is a
purely transcendental extension of infinite transcendence degree
over the rational function field $F _{n-1}(T _{n})$.
\par
{\rm (b)} The maximal separable (algebraic) extension $E _{n}$ of
$\Lambda _{n}$ in $F _{n}$ satisfies the equalities $[E _{n}\colon
E _{n} ^{q}] = \infty $, {\rm Brd}$_{p}(E _{n}) = {\rm abrd}_{p}(E
_{n}) = [n/2]$, for all
\par\noindent
$p \in \mathbb{P}_{F _{0}}$, and {\rm Brd}$_{q}(E _{n} ^{\prime })
= n - 1$, for every finite extension $E _{n} ^{\prime }/E _{n}$.
\end{prop}
\par
\smallskip
\begin{proof}
The assertion of Proposition \ref{prop4.9} (a) is known (cf.,
e.g., \cite{BlKu}), and it implies $[E _{n}\colon E _{n} ^{q}] =
\infty $. Let $w _{n}$ be the natural discrete valuation of $F
_{n}$ trivial on $F _{n-1}$, and $v _{n}$ be the valuation of $E
_{n}$ induced by $w _{n}$. Then $(F _{n}, w _{n})$ is complete and
$E _{n}$ is dense in $F _{n}$, whence, $v _{n}(E _{n}) = w _{n}(F
_{n})$ and $F _{n-1}$ is the residue field of $(E _{n}, v _{n})$
and $(F _{n}, w _{n})$ (cf. \cite[Ch.~XII, \S{5}]{L}); in
particular, $v _{n}$ is discrete. Similarly, if $n \ge 2$, then
the natural $\mathbb{Z} ^{n}$-valued valuation $\theta _{n} '$ of
$F _{n}$ (trivial on $F _{0}$) is Henselian and induces on $E
_{n}$ a valuation $\theta _{n}$. Also, $v _{n}$ is Henselian (cf.
\cite[Corollary~18.3.3]{Ef}), and $\theta _{n}$ extends the
natural $\mathbb{Z} ^{n-1}$-valued valuation $\theta _{n-1}'$ of
$F _{n-1}$. As $\theta _{n-1}'$ is Henselian and $F _{n-1}(T _{n})
\subset E _{n}$, this ensures that so is $\theta _{n}$ (see
\cite[Proposition~A.15]{TW}), $F _{0}$ is the residue field of $(E
_{n}, \theta _{n})$, and $\theta _{n}(E _{n}) = \mathbb{Z} ^{n}$.
At the same time, it follows from Cohn's theorem (cf.
\cite[Theorem~1]{Cohn} or \cite[Proposition~1.16]{TW}) and the
Henselian property of $v _{n}$ that Brd$_{p}(E _{n}) \le {\rm
Brd}_{p}(F _{n})$, for each $p$. In addition, $(F _{n}, \theta
_{n}')$ is maximally complete with a residue field $F _{0}$ (cf.
\cite[Theorem~18.4.1]{Ef}), whence, by \cite[Proposition~3.5]{Ch6},
Brd$_{q}(F _{n}) = {\rm abrd}_{q}(F _{n}) = n - 1$. Since
\par\vskip0.063truecm\noindent
$\theta _{n}(E _{n} ^{\prime }) \cong \mathbb{Z} ^{n} \cong \theta
_{n}'(F _{n} ^{\prime })$ and $F _{0}$ is the residue field of
both $(E _{n} ^{\prime }, \theta _{n,E _{n}'})$ and $(F _{n},
\theta _{n,F _{n}'}')$ whenever $E _{n} ^{\prime }/E _{n}$ and $F
_{n} ^{\prime }/F _{n}$ are finite extensions, one obtains from
\cite[Proposition~5.3 (b)]{Ch7} and \cite[Lemma~4.2]{Ch5} that
Brd$_{q}(E _{n} ^{\prime }) \ge n - 1$
\par\vskip0.063truecm\noindent
and Brd$_{p}(E _{n} ^{\prime }) = {\rm Brd}_{p}(F _{n} ^{\prime })
= [n/2]$, for each $p \in \mathbb{P}_{F_{0}}$. Note finally that
$v _{n}(E _{n} ^{\prime }) \cong \mathbb{Z} \cong w _{n}(F _{n}
^{\prime })$, $v _{n,E_{n}}$ is Henselian, and the completion of
$E _{n} ^{\prime }$ with respect to the topology of $v _{n,E
_{n}'}$ is a finite extension of $F _{n}$ (cf. \cite[Ch. XII,
Proposition~3.1]{L}), so Cohn's theorem implies Brd$_{q}(E _{n}
^{\prime }) = n - 1$, $n \in \mathbb{N}$.
\end{proof}
\par
\medskip
\section{\bf Lemmas on $p$-powers and finite-dimensional central
subalgebras of division {\rm LBD}-algebras}
\par
\medskip
Let $R$ be a central division LBD-algebra over a virtually perfect
field $K$ with abrd$_{p}(K) < \infty $, $p \in \mathbb{P}$. The
existence of finite $p$-powers $k(p)$ of $R/K$, $p \in \mathbb{P}$
(by \cite[Lemma~3.9]{Ch2}), imposes essential restrictions on a
number of algebraic properties of $R$, especially, on the
extensions of $K$ which embed in $R$ as $K$-subalgebras. For
example, it turns out that if $K(p) \neq K$, for some $p > 2$,
then $K(p)/K$ is an infinite extension (the additive group
$\mathbb{Z} _{p}$ of $p$-adic integers, endowed with its natural
topology, is a homomorphic image of $\mathcal{G}(K(p)/K)$, see
\cite{Wh}), whence, $K(p)$ is not isomorphic to a $K$-subalgebra
of $R$. In this Section we present results on $p$-powers and
$p$-splitting fields, obtained in the case of dim$(K _{\rm sol})
\le 1$. These results provide the basis of the proofs of Theorems
\ref{theo3.1}, \ref{theo3.2} and Lemma \ref{lemm3.3}. The first
one is an immediate consequence of \cite[Lemmas~3.12 and
3.13]{Ch2} and can be stated as follows:
\par
\medskip
\begin{lemm}
\label{lemm5.1} Assume that $R$ is a central division {\rm
LBD}-algebra over a virtually perfect field $K$ with {\rm dim}$(K
_{\rm sol}) \le 1$ and {\rm abrd}$_{p}(K) < \infty $, $p \in
\mathbb{P}$. Let $K ^{\prime }/K$ be a finite extension, $R
^{\prime }$ the underlying (central) division $K ^{\prime
}$-algebra of the {\rm LBD}-algebra $R \otimes _{K} K ^{\prime }$,
$\gamma $ the integer for which $R \otimes _{K} K ^{\prime }$ and
the matrix ring $M _{\gamma }(R ^{\prime })$ are isomorphic as $K
^{\prime }$-algebras, and for each $p \in \mathbb{P}$, let $k(p)$
and $k(p)'$ be the $p$-powers of $R/K$ and $R ^{\prime }/K
^{\prime }$, respectively. Then:
\par
{\rm (a)} The greatest integer $\mu (p) \ge 0$ for which $p ^{\mu
(p)} \mid \gamma $ is equal to $k(p) - k(p)'$; hence, $k(p) \ge
k(p)'$ and $p ^{1+k(p)} \nmid \gamma $, for any $p \in \mathbb{P}$;
\par
{\rm (b)} The equality $k(p) = k(p)'$ holds if and only if $p
\nmid \gamma $; specifically, if $k(p) = 0$, then $k(p)'= 0$ and
$p \nmid \gamma $.
\par
{\rm (c)} $K ^{\prime }$ is a $p$-splitting field of $R/K$ if and
only if $k(p)' = 0$, that is, $p \nmid [K ^{\prime }(r')\colon K
^{\prime }]$, for any $r' \in R ^{\prime }$.
\end{lemm}
\par
\medskip
As a matter of fact, Lemma \ref{lemm5.1} (a) is identical in
content with \cite[Lemmas~3.12 and 3.13]{Ch2}, and it implies
Lemma \ref{lemm5.1} (b) and (c).
\par
\medskip
\begin{rema}
\label{rema5.2} The proofs of \cite[Lemmas~3.12 and 3.13]{Ch2}
rely on the condition that dim$(K _{\rm sol}) \le 1$, more
precisely, on its restatement that abrd$_{p}(K _{p}) = 0$, for
each $p \in \mathbb{P}$, where $K _{p}$ is the fixed field of a
Hall pro-$(\mathbb{P} \setminus \{p\})$-subgroup $H _{p}$ of
$\mathcal{G}(K _{\rm sol}/K)$. It is not known whether the
assertions of Lemma \ref{lemm5.1} remain valid if this condition
is dropped. Also, Koenigsmann's question of whether {\rm dim}$(E
_{\rm sol}) \le 1$, for every field $E$ (posed in \cite{Koe}) is
open; an affirmative answer is predicted by the
Bogomolov-Positselski conjecture (see \cite[Conjecture~1]{Po}).
Here we note that the conclusion of Lemma \ref{lemm5.1} holds if
the assumption that dim$(K _{\rm sol}) \le 1$ is replaced by the
one that $R \otimes _{K} K ^{\prime }$ is a division $K ^{\prime
}$-algebra; then it follows from \cite[Proposition~3.3]{Ch2} that
$k(p) = k(p)'$, for every $p \in \mathbb{P}$.
\end{rema}
\par
\smallskip
\begin{lemm}
\label{lemm5.3} Assuming that $K$ and $R$ satisfy the conditions
of Lemma \ref{lemm5.1}, let $D \in d(K)$ be a $K$-subalgebra of
$R$, and for each $p \in \mathbb{P}$, let $k(p)$ and $k(p)'$ be
the $p$-powers of $R/K$ and $C _{R}(D)/K$, respectively. Then:
\par
{\rm (a)} The sequence $k(p) - k(p)'\colon p \in \mathbb{P}$,
consists of integers $\ge 0$, which are equal to the power of $p$
in the primary decomposition of {\rm deg}$(D)$, for each $p$;
\par
{\rm (b)} $k(p) = k(p)'$ if and only if $p \nmid {\rm deg}(D)$; in
this case, a finite extension $K ^{\prime }$ of $K$ is a
$p$-splitting field of $R/K$ if and only if so is $K ^{\prime }$ for
$C _{R}(D)/K$;
\par
{\rm (c)} If $k(p)' = 0$, for some $p \in \mathbb{P}$, then a
finite extension $K ^{\prime }$ of $K$ is a $p$-splitting field of
$R/K$ if and only if $p \nmid {\rm ind}(D \otimes _{K} K ^{\prime
})$.
\end{lemm}
\par
\smallskip
\begin{proof}
It is known (cf. \cite[Sect. 13.1, Corollary~b]{P}) that if $K
_{1}$ is a maximal subfield of $D$, then $K _{1}/K$ is a field
extension, $[K _{1}\colon K] = {\rm deg}(D) := d$ and $D \otimes
_{K} K _{1} \cong M _{d}(K _{1})$ as $K _{1}$-subalgebras. Also,
by the double centralizer theorem (see \cite[Theorems~4.3.2 and
4.4.2]{He}), $R = D \otimes _{K} C _{R}(D)$ and $C _{R}(D) \otimes
_{K} K _{1}$ is a central division $K _{1}$-algebra equal to $C
_{R}(K _{1})$. In view of \cite[Propositions~3.1 and 3.3]{Ch2},
this ensures that $k(p)'$ equals the $p$-power of $(C _{R}(D)
\otimes _{K} K _{1})/K _{1}$, for each $p \in \mathbb{P}$.
Applying now Lemma \ref{lemm5.1}, one proves Lemma \ref{lemm5.3}
(a). Lemma \ref{lemm5.3} (b)-(c) follows from Lemmas \ref{lemm5.1}
and \ref{lemm5.3} (a), combined with \cite[Lemma~3.5]{Ch2} and
\cite[Sect. 9.3, Corollary~b]{P}.
\end{proof}
\par
\medskip
The following lemma (for a proof, see \cite[Lemma~7.4]{Ch2}) can
be viewed as a generalization of the uniqueness part of the
primary tensor product decomposition theorem for algebras $D \in
d(K)$ over an arbitrary field $K$.
\par
\medskip
\begin{lemm}
\label{lemm5.4} Let $\Pi $ be a finite subset of $\mathbb P$, and
let $S _{1}$, $S _{2}$ be central division {\rm LBD}-algebras over
a field $K$ with {\rm abrd}$_{p}(K) < \infty $, for all $p \in
\mathbb P$. Assume that $k(p)_{S _{1}} = k(p)_{S _{2}} = 0$, $p
\in \Pi $, the $K$-algebras $A _{1} \otimes _{K} S _{1}$ and $A
_{2} \otimes S _{2}$ are $K$-isomorphic, where $A _{i} \in s(K)$,
$i = 1, 2$, and {\rm deg}$(A _{1}){\rm deg}(A _{2})$ is not
divisible by any $\bar p \in \mathbb P \setminus \Pi $. Then $A
_{1} \cong A _{2}$ as $K$-algebras.
\end{lemm}
\par
\medskip
For a proof of our next lemma, we refer the reader to
\cite[Lemmas~8.3 and 8.4]{Ch2}, which have been proved under the
assumption that $R$ is a central division LBD-algebra over a field
$K$ of arithmetic type. Therefore, we note that the proof in \cite{Ch2}
remains valid if the assumption on $K$ is replaced by the one that
abrd$_{p}(K) < \infty $, $p \in \mathbb{P}$, dim$(K _{\rm sol})
\le 1$, $K$ is virtually perfect, and there exist $p$-splitting
fields $E _{p}: p \in \mathbb{P}$, of $R/K$ with $E _{p} \subseteq
K(p)$, for each $p$.
\par
\medskip
\begin{lemm}
\label{lemm5.5} Let $R$ be a central division {\rm LBD}-algebra
over a field $K$ with {\rm dim}$(K _{\rm sol}) \le 1$, and for
each $p \in \mathbb{P}$, let $k(p)$ be the $p$-power of $R/K$, and
$E _{p}$ be a $p$-splitting field of $R/K$, such that $E _{p} \in
I(K(p)/K)$. Then:
\par
{\rm (a)} The matrix ring $M _{\gamma (p)}(R)$, where $\gamma (p)
= [E _{p}\colon K].p ^{-k(p)}$, is an Artinian central simple {\rm
LBD}-algebra over $K$, which possesses a subalgebra $\Delta _{p}
\in s(K)$, such that {\rm deg}$(\Delta _{p}) = [E _{p}\colon K]$
and $E _{p}$ is embeddable in $\Delta _{p}$ as a $K$-subalgebra.
Moreover, if $[E _{p}\colon K] = p ^{k(p)}$, i.e. $E _{p}$ is
embeddable in $R$ as a $K$-subalgebra, then $\Delta _{p}$ is a
$K$-subalgebra of $R$.
\par
{\rm (b)} The centralizer of $\Delta _{p}$ in $M _{\gamma (p)}(R)$ is
a central division $K$-algebra of $p$-power zero.
\end{lemm}
\par
\medskip
The following lemma generalizes \cite[Lemma~8.5]{Ch2} to the case
where $K$ and $R$ satisfy the conditions of Lemma \ref{lemm5.5}.
For this reason, we take into account that the proof of the lemma
referred to, given in \cite{Ch2}, remains valid under the noted
weaker conditions. Our next lemma can also be viewed as a
generalization of the well-known fact that, for any field $E$, $D
_{1} \otimes _{E} D _{2} \in d(E)$ whenever $D _{i} \in d(E)$, $i
= 1, 2$, and $\gcd \{{\rm deg}(D _{1}), {\rm deg}(D _{2})\} = 1$
(see \cite[Sect. 13.4]{P}). Using this lemma and the uniqueness
part of the Wedderburn-Artin theorem, one obtains that, in the
setting of Lemma \ref{lemm5.5}, the underlying central division
$K$-algebra of $\Delta _{p}$ is embeddable in $R$ as a
$K$-subalgebra.
\par
\medskip
\begin{lemm}
\label{lemm5.6} Let $K$ be a field, $R$ a central division {\rm
LBD}-algebra over $K$, and $E _{p}$, $p \in \mathbb{P}$, be
extensions of $K$ satisfying the conditions of Lemma
\ref{lemm5.5}. Also, let $D \in d(K)$ be a division $K$-algebra
such that $\gcd \{{\rm deg}(D), [K(\alpha )\colon K]\} = 1$, for
each $\alpha \in R$. Then $D \otimes _{K} R$ is a central division
{\rm LBD}-algebra over $K$.
\end{lemm}
\par
\medskip
The following two lemmas give us the possibility to deduce Lemma
\ref{lemm3.3} from Lemma \ref{lemm2.2} by the method of proving
\cite[Lemma~8.3]{Ch2} (see (8.6)).
\par
\medskip
\begin{lemm}
\label{lemm5.7} Let $K$, $R$ and $k(p)\colon p \in \mathbb{P}$,
satisfy the conditions of Lemma \ref{lemm5.1}, and let $K _{1}$,
$K _{2}$ and $E$ be finite extensions of $K$ in an algebraic
closure of $K _{\rm sep}$. Denote by $R _{1}$ and $R _{2}$ the
underlying division algebras of $R \otimes _{K} K _{1}$ and $R
\otimes _{K} K _{2}$, respectively, and suppose that there exist
algebras $D _{1} \in d(K _{1})$, $D _{2} \in d(K _{2})$, such that
$D _{i}$ is a $K _{i}$-subalgebra of $R _{i}$ and {\rm deg}$(D
_{i}) = p _{0} ^{k(p_{0})}$, for a given $p _{0} \in \mathbb{P}$
and each index $i$. Then:
\par
{\rm (a)} The underlying division $K _{1}K _{2}$-algebras of
$R _{1} \otimes _{K _{1}} K _{1}K _{2}$, $R _{2} \otimes _{K _{2}} K
_{1}K _{2}$ and $R \otimes _{K} K _{1}K _{2}$ are isomorphic;
\par
{\rm (b)} $p _{0} \nmid [K _{i}(c _{i})\colon K _{i}]$, for any $c
_{i} \in C _{R _{i}}(D _{i})$, and $i = 1, 2$.
\par
{\rm (c)} If $p _{0} \nmid [K _{1}K _{2}\colon K]$, then $D _{1}
\otimes _{K _{1}} K _{1}K _{2}$ and $D _{2} \otimes _{K _{2}} K
_{1}K _{2}$ are isomorphic central division $K _{1}K
_{2}$-algebras; for example, this holds in case $p _{0} \nmid [K
_{i}\colon K]$, $i = 1, 2$, and $\gcd \{[K _{1}\colon K _{0}], [K
_{2}\colon K _{0}]\} = 1$, where $K _{0} = K _{1} \cap K _{2}$.
\end{lemm}
\par
\smallskip
\begin{proof}
The $K _{1}K _{2}$-algebras $R \otimes _{K} K _{1}K _{2}$ and $(R
\otimes _{K} K _{i}) \otimes _{K _{1}} K _{1}K _{2}$, $i = 1, 2$,
are isomorphic, by known general properties of tensor products
(cf. \cite[Sect.~9.4, Corollary~a]{P}). Since the considered
algebras are central simple and Artinian, this fact enables one to
deduce Lemma \ref{lemm5.7} (a) from Wedderburn-Artin's theorem and
\cite[Sect. 9.3, Corollary~b]{P}. In addition, it follows from
Lemma \ref{lemm5.1}, the assumptions on $D _{1}$ and $D _{2}$, and
the double centralizer theorem that $k(p)$ equals the $p$-powers
of $R _{i}/K _{i}$, and $C _{R _{i}}(D _{i})$ is a central
division $K _{i}$-subalgebra of $R _{i}$, for each $i$. Hence, by
Lemma \ref{lemm5.3} (a), $C _{R _{1}}(D _{1})/K _{1}$ and $C _{R
_{2}}(D _{2})/K _{2}$ are of $p _{0}$-power zero, which proves
Lemma \ref{lemm5.7} (b).
\par
We turn to the proof of Lemma \ref{lemm5.7} (c). Assume that $p
_{0} \nmid [K _{1}K _{2}\colon K]$ and denote by $R ^{\prime }$
the underlying division $K _{1}K _{2}$-algebra of $R \otimes _{K}
K _{1}K _{2}$. Then, by Lemma \ref{lemm5.1}, $k(p _{0})$ equals
the $p _{0}$-power of $R ^{\prime }/K _{1}K _{2}$. Applying
\cite[Lemma~3.5]{Ch2} (or results of \cite[Sect. 13.4]{P}), one
also obtains that $D _{i} \otimes _{K _{i}} K _{1}K _{2} \in d(K
_{1}K _{2})$ and $D _{i} \otimes _{K _{i}} K _{1}K _{2}$ is
isomorphic to a $K _{1}K _{2}$-subalgebra $D _{i} ^{\prime }$ of
$R ^{\prime }$, for $i = 1, 2$.
\par\vskip0.04truecm\noindent
As above, now it follows that, for each $i$, $R ^{\prime } = D
_{i} ^{\prime } \otimes _{K _{1}K _{2}} C _{R'}(D _{i} ^{\prime
})$, $C _{R'}(D _{i} ^{\prime })$
\par\vskip0.04truecm\noindent
is a central division $K _{1}K _{2}$-algebra, and $C _{R'}(D _{i}
^{\prime })/K _{1}K _{2}$ is of zero $p _{0}$-power;
\par\vskip0.04truecm\noindent
thus $p _{0} \nmid [K _{1}K _{2}(c')\colon K _{1}K _{2}]$, for any
$c' \in C _{R'}(D _{1} ^{\prime }) \cup C _{R'}(D _{2} ^{\prime
})$. Therefore, by
\par\vskip0.04truecm\noindent
Lemma \ref{lemm5.4}, $D _{1} ^{\prime } \cong D _{2} ^{\prime }$,
whence, $D _{1} \otimes _{K _{1}} K _{1}K _{2} \cong D _{2}
\otimes _{K _{2}} K _{1}K _{2}$ over $K _{1}K _{2}$. The latter
part of Lemma \ref{lemm5.7} (c) is obvious, so our proof is
complete.
\end{proof}
\par
\medskip
\begin{lemm}
\label{lemm7.7} Assume that $E$ is a field and $B$, $W$, $W
^{\prime }$ are finite extensions of $E$ in $E _{\rm sep}$, such
that $W \in I(W ^{\prime }/E)$ and $W ^{\prime }/E$ is a Galois
extension with $\mathcal{G}(W ^{\prime }/E)$ nilpotent. Then $BW$
is an extension of $W$ and $[BW\colon W] \mid [B\colon E]$.
\end{lemm}
\par
\smallskip
\begin{proof}
It is easily verified that $BW/E$ is a finite extension with
\par\vskip0.032truecm\noindent
$[BW\colon E] = [BW\colon W][W\colon E] = [BW\colon B][B\colon
E]$. In addition, it follows
\par\vskip0.032truecm\noindent
from Galois theory (see \cite[Ch. VI, Theorem~1.12]{L}) that if
$W/E$ is a Galois extension, then $BW/B$ is Galois and $[BW\colon
B] \mid [W\colon E]$, which implies $[BW\colon W] \mid [B\colon
E]$, as required. Henceforth, we suppose that $W/E$ is not Galois;
in particular, this yields $W \neq E$ and $\mathcal{G}(W ^{\prime
}/W)$ is a proper subgroup of $\mathcal{G}(W ^{\prime }/E)$. Note
further that, by the Burnside-Wielandt theorem (see
\cite[Theorem~17.1.4]{KM}) and the nilpotency of $\mathcal{G}(W
^{\prime }/E)$, $\mathcal{G}(W ^{\prime }/W)$ is subnormal in
$\mathcal{G}(W/E)$, which allows to deduce from Galois theory that
$E$ has a cyclic extension $W _{0}$ in $W$ of prime degree $[W
_{0}\colon E]$. Therefore, by the preceeding part of our proof,
$BW _{0}/W _{0}$ is a field extension, $[BW _{0}\colon W _{0}]
\mid [B\colon E]$ and $\mathcal{G}(W ^{\prime }/W _{0})$ is a
maximal subgroup of $\mathcal{G}(W ^{\prime }/E)$ of index $[W
_{0}\colon E]$. It is now easy to see that the field $W _{0}$ and
its extensions $BW _{0}, W$ and $W ^{\prime }$ satisfy the
conditions of Lemma \ref{lemm7.7}, so the assertion of this lemma
requires that $[BW\colon W] \mid [BW _{0}\colon W _{0}]$. Since
$[W ^{\prime }\colon W _{0}] < [W ^{\prime }\colon E]$ and $[BW
_{0}\colon W _{0}] \mid [B\colon E]$, these observations make it
possible to complete the proof of Lemma \ref{lemm7.7}, by assuming
the opposite; one arrives at a contradiction by taking as a
counter-example fields $E$, $B$, $W$ and $W ^{\prime }$, chosen so
that $[W ^{\prime }\colon E]$ be minimal.
\end{proof}
\par
\medskip
As shown at the end of this section, the next lemma plays a major
role in the reduction of Theorems \ref{theo3.1} and \ref{theo3.2}
to consequences of Lemma \ref{lemm3.3}.
\par
\medskip
\begin{lemm}
\label{lemm5.9} Let $(K, v)$ be a Henselian field with {\rm
dim}$(K _{\rm sol}) \le 1$ and {\rm abrd}$_{\ell }(K)$ finite, for
all $\ell \in \mathbb{P}$. Fix $p \in \mathbb{P}$ and a field $M
\in I(M ^{\prime }/K)$, for some finite Galois extension $M
^{\prime }$ of $K$ in $K _{\rm sep}$ with $\mathcal{G}(M ^{\prime
}/K)$ nilpotent and $[M ^{\prime }\colon K]$ not divisible by $p$.
Assume that $R$ is a central division {\rm LBD}-algebra over $K$,
$R _{M}$ is the underlying division $M$-algebra of $R \otimes _{K}
M$, and $R _{M}$ has an $M$-subalgebra $\Delta _{M}$, such that
the following (equivalent) conditions hold:
\par
{\rm (c)} $M$ is a $p'$-splitting field of $R/K$, for every $p'
\in \mathbb{P}$ dividing $[M\colon K]$; $\Delta _{M} \in d(M)$ and
{\rm deg}$(\Delta _{M}) = p ^{k(p)}$, where $k(p)$ is the
$p$-power of $R/K$;
\par
{\rm (cc)} $\gcd \{p[M\colon K], [M(z _{M})\colon M]\} = 1$, for
every $z _{M} \in C _{R _{M}}(\Delta _{M})$.
\par\smallskip\noindent
Then $\Delta _{M} \cong \Delta \otimes _{K} M$ as $M$-algebras,
for some subalgebra $\Delta \in d(K)$ of $R$.
\end{lemm}
\par
\smallskip
\begin{proof}
The equivalence of conditions (c) and (cc) follows from Lemmas
\ref{lemm5.1} and \ref{lemm5.3}. Note further that if $M \neq K$,
then $M$ contains as a subfield a cyclic extension $M _{0}$ of $K$
of degree $p' \neq p$. Since, by the Burnside-Wielandt theorem,
maximal subgroups of nilpotent finite groups are normal, this can
be deduced from Galois theory. Considering $M _{0}$ and the
underlying division $M _{0}$-algebra $R _{0}$ of $R \otimes _{K} M
_{0}$, instead of $K$ and $R$, respectively, and taking into
account that $R \otimes _{K} M \cong (R \otimes _{K} M _{0})
\otimes _{M _{0}} M$ as $M$-algebras, one concludes that
conditions (c) and (cc) of Lemma \ref{lemm5.9} are fulfilled
again. Therefore, a standard inductive argument shows that it
suffices to prove Lemma \ref{lemm5.9} under the extra hypothesis
that there exists a subalgebra $\Delta _{0} \in d(M _{0})$ of $R
_{0}$, such that $\Delta _{0} \otimes _{M _{0}} M \cong  \Delta
_{M}$ as $M$-algebras. Let $\varphi $ be a $K$-automorphism of $M
_{0}$ of order $p'$, and let $\bar \varphi $ be the unique
$K$-automorphism of $R \otimes _{K} M _{0}$ extending $\varphi $
and acting on $R$ as the identity. Then it follows from the
Skolem-Noether theorem (cf. \cite[Theorem~4.3.1]{He}) and from the
existence of an $M _{0}$-isomorphism $R \otimes _{K} M _{0} \cong
M _{p*}(M _{0}) \otimes _{M _{0}} R _{0}$ (where $p* = p$ or $p* =
1$
\par\vskip0.04truecm\noindent
depending on whether or not $M _{0}$ embeds in $R$ as a
$K$-subalgebra) that $R _{0}$ has a $K$-automorphism $\tilde
\varphi $ extending $\varphi $. Note also that $p \nmid [M _{0}(z
_{0})\colon M _{0}]$, for any $z _{0} \in C _{R _{0}}(\Delta
_{0})$. This is implied by Lemma \ref{lemm5.1}, condition (cc) of
Lemma \ref{lemm5.9}, and the fact that $C _{R _{M}}(\Delta _{M})$
is the underlying division $M$-algebra of $C _{R _{M _{0}}}(\Delta
_{0}) \otimes _{M _{0}} M$. Hence, by Lemma \ref{lemm5.4}, $\Delta
_{0}$ is $M _{0}$-isomorphic to its image $\Delta _{0} ^{\prime }$
under $\tilde \varphi $, so it follows from Skolem-Noether's
theorem that $\varphi $ extends to a $K$-automorphism of $\Delta
_{0}$. As deg$(\Delta _{0}) = {\rm deg}(\Delta ) = p ^{k(p)}$ and
$p \nmid [M _{0}\colon K]$, this enables one to deduce from
Teichm\"{u}ller's theorem (cf. \cite[Sect.~9, Theorem~4]{Dr1}) and
\cite[Lemma~3.5]{Ch2} the existence of an $M _{0}$-isomorphism
\par\vskip0.04truecm\noindent
$\Delta _{0} \cong \Delta \otimes _{K} M _{0}$, for some
$K$-subalgebra $\Delta  \in d(K)$ of $R$.
\end{proof}
\par
\smallskip
We are now prepared to deduce Theorems \ref{theo3.1} and
\ref{theo3.2} from Lemma \ref{lemm3.3}. Let $K$ be a virtually
perfect field with dim$(K _{\rm sol}) \le 1$ and abrd$_{\ell }(K)
< \infty $, for each $\ell \in \mathbb{P}_{K}$, and let $R$ be a
central division LBD-algebra over $K$. Denote by $k(p)$ the
$p$-power of $R/K$, for each $p \in \mathbb{P}$, and suppose that
$K$ has finite extensions $E _{p} \in I(K(p)/K)\colon p \in
\mathbb{P}$, which are $p$-splitting fields of $R/K$. Then one
obtains from Lemmas \ref{lemm5.5}, \ref{lemm5.6} and
\ref{lemm7.7}, by the method of proving \cite[Lemma~8.3]{Ch2},
that for each $p \in \mathbb{P}$, there exists a finite subset
$\Pi _{p}$ of $\mathbb{P} \setminus \{p\}$, such that the
compositum $K _{\Pi _{p}}$ of the fields $E _{p'}\colon p' \in \Pi
_{p}$, and the underlying division $K _{\Pi _{p}}$-algebra
$\mathcal{R}_{\Pi _{p}}$ of $R \otimes _{K} K _{\Pi _{p}}$, have
the following properties: $\mathcal{R}_{\Pi _{p}}$ possesses a
subalgebra $\Delta _{K _{\Pi _{p}}} \in d(K _{\Pi _{p}})$ of
degree $p ^{k(p)}$; $E _{p}K _{\Pi _{p}}$ is a splitting field of
$\Delta _{K_{\Pi _{p}}}$ and a $p$-splitting field of
$\mathcal{R}_{\Pi _{p}}$. In view of Lemma \ref{lemm5.9}, this
means that $R$ has a unique, up-to isomorphism, $K$-subalgebra $R
_{p} \in d(K)$ of degree $p ^{k(p)}$; also, it follows from
\cite[Lemma~3.5]{Ch2} that $E _{p}$ is a spliting field of $R
_{p}$. Applying now Lemma~3.2 of \cite{Ch9} (which is in fact a
special case of \cite[Lemma~3.5]{Ch2}), one obtains that $R$
contains as a $K$-subalgebra an isomorphic copy of $\otimes _{p
\in \Pi } R _{p}$, for each finite subset $\Pi $ of $\mathbb{P}$,
where $\otimes = \otimes _{K}$. Now it follows from Lemma
\ref{lemm5.3} and Skolem-Noether's theorem (cf.
\cite[Theorem~4.3.1]{He}) that the algebras $R _{p}$, $p \in
\mathbb{P}$, can be chosen so that $R _{p'} \subseteq C _{R}(R
_{p''})$ whenever $p', p'' \in \mathbb{P}$ and $p' \neq p''$.
Therefore, $R$ has $K$-subalgebras $T _{n}$, $n \in \mathbb{N}$,
such that $T _{n} \cong \otimes _{j=1} ^{n} R _{p _{j}}$ and $T
_{n} \subseteq T _{n+1}$, for each $n$; here $\otimes = \otimes
_{K}$ and $\mathbb{P}$ is presented as a growing sequence $p
_{n}\colon n \in \mathbb{N}$. Clearly, the union $\widetilde R =
\cup _{n=1} ^{\infty } T _{n} := \otimes _{n=1} ^{\infty } R _{p
_{n}}$ is a central $K$-subalgebra of $R$. Note further that $R =
T _{n} \otimes _{K} C _{R}(T _{n})$, for every $n \in \mathbb{N}$
(apply \cite[Theorem~4.4.2]{He}), which enables one to deduce from
Lemmas \ref{lemm5.1}, \ref{lemm5.4}, \ref{lemm5.6}, and
\cite[Lemma~3.5]{Ch2}, that a $K$-subalgebra $T$ of $R$ with
$[T\colon K] < \infty $ embeds in the $K$-algebra $T _{n}$ if $p
_{n'} \nmid [T\colon K]$, for any $n' > n$; in
\par\vskip0.04truecm\noindent
particular, $T$ embeds in $\widetilde R$. One also sees that $K =
\cap _{n=1} ^{\infty } C _{R}(T _{n}) = C _{R}(\widetilde R)$.
\par
\medskip
\section{\bf Henselian fields $(K, v)$ with char$(\widehat K) = q >
0$ and {\rm abrd}$_{q}(K(q)) \le 1$}
\par
\medskip
The question of whether abrd$_{q}(\Phi (q)) = 0$, for every field
$\Phi $ of characteristic $q > 0$ seems to be open. This Section
gives a criterion for a Henselian field $(K, v)$ with
char$(\widehat K) = q$ and $\widehat K$ of arithmetic type to
satisfy the equality abrd$_{q}(K(q)) = 0$. To prove this criterion
we need the following two lemmas.
\par
\medskip
\begin{lemm}
\label{lemm6.1} Let $(K, v)$ be a Henselian field with {\rm
char}$(\widehat K) = q > 0$ and
\par\noindent
$\widehat K \neq \widehat K ^{q}$, and in case {\rm char}$(K) =
0$, suppose that $v$ is discrete and
\par\noindent
$v(q) \in qv(K)$. Let also
$\widetilde \Lambda /\widehat K$ be an inseparable extension
of degree $q$. Then there is $\Lambda \in I(K(q)/K)$, such that
$[\Lambda \colon K] = q$ and $\widehat \Lambda $ is $\widehat
K$-isomorphic to $\widetilde \Lambda $.
\end{lemm}
\par
\smallskip
\begin{proof}
When char$(K) = 0$, our assertion is contained in
\cite[Lemma~5.4]{Ch8}, so we assume here that char$(K) = q$. It
follows from the condition on $\widetilde \Lambda /\widehat K$
that there exists $\tilde a \in \widehat K \setminus \widehat K
^{q}$, such that $\widetilde \Lambda = \widehat K(\sqrt[q]{\tilde
a})$. Hence, by the Artin-Schreier theorem (cf. \cite[Ch. VI,
\S{6}]{L}), if $\pi \in K$, $v(\pi ) > 0$ and $a \in O _{v}(K)$ is
chosen so that $\hat a = \tilde a$, then one may take as $\Lambda
$ the extension of $K$ in $K _{\rm sep}$ obtained by adjunction of
a root of the polynomial $X ^{q} - X - a\pi ^{-q}$.
\end{proof}
\par
\medskip
\begin{lemm}
\label{lemm6.2} Let $(K, v)$ be a Henselian field, $L/K$ an
inertial extension, and $N(L/K)$ the norm group of $L/K$. Then
$\nabla _{0}(K)$ is a subgroup of $N(L/K)$.
\end{lemm}
\par
\smallskip
\begin{proof}
This is a special case of \cite[Proposition~2]{Er}.
\end{proof}
\par
\medskip
Next we show that a Henselian field $(K, v)$ with char$(\widehat
K) = q > 0$ satisfies abrd$_{q}(K(q)) = 0$, provided that char$(K)
= q$ or the valuation $v$ is discrete. In view of the
Albert-Hochschild theorem, the equality abrd$_{q}(K(q)) = 0$ ensures
that Br$(K(q) ^{\prime }) _{q} = \{0\}$, for every finite extension
$K(q) ^{\prime }/K(q)$.
\par
\medskip
\begin{lemm}
\label{lemm6.3} Let $(K, v)$ be a Henselian field with {\rm
char}$(\widehat K) = q > 0$, and in case {\rm char}$(K) = 0$, let
$v$ be discrete. Then $v(K(q)) = qv(K(q))$, the residue field
$\widehat {K(q)} = \widehat K(q) _{\rm ins}$ of $(K(q), v
_{K(q)})$ is perfect, and {\rm abrd}$_{q}(K(q)) = 0$.
\end{lemm}
\par
\medskip
Since the proof of Lemma \ref{lemm6.3} relies on the
presentability of cyclic \par\noindent $K$-algebras of degree $q$
as $q$-symbol algebras over $K$, we recall some basic facts about
$q$-symbol algebras over any field $E$ with char$(E) = q > 0$.
Firstly, for each pair $a \in E$, $b \in E ^{\ast }$, we have $[a,
b) _{E} \in s(E)$ and deg$([a, b) _{E}) = q$ (cf.
\cite[Corollary~2.5.5]{GiSz}). Secondly, if $[a, b) _{E} \in
d(E)$, then the polynomial
\par\noindent
$f _{a}(X) = X ^{q} - X - a \in E[X]$ is irreducible over $E$.
This follows from Artin-Schreier's theorem which also shows that
if $f _{a}(X)$ is irreducible over $E$ and $\xi \in E _{\rm sep}$
is a root of $f _{a}(X)$, then $E(\xi )/E$ is a degree $q$ cyclic
field extension and $[a, b) _{E}$ is isomorphic to the cyclic
$E$-algebra $(E(\xi )/E, \sigma , b)$, where $\sigma $ is the
$E$-automorphism of $E(\xi )$ mapping $\xi $ into $\xi + 1$;
hence, by \cite[Sect. 15.1, Proposition~b]{P}, $[a, b) _{E} \in
d(E)$ if and only if $b \notin N(E(\xi )/E)$.
\par
\medskip\noindent
{\it Proof of Lemma \ref{lemm6.3}.} It is clear from Galois
theory, the definition of $K(q)$ and the closeness of the class of
pro-$q$-groups under the formation of profinite group extensions
that $\widetilde K(q) = K(q)$, for every $\widetilde K \in
I(K(q)/K)$; in particular, $K(q)(q) = K(q)$, which means that
$K(q)$ does not admit cyclic extensions of degree $q$. As $(K(q),
v_{K(q)})$ is Henselian, this allows to deduce from
\cite[Lemma~4.2]{Ch5} and \cite[Lemma~2.2]{Ch8} that $v(K(q)) =
qv(K(q))$. We show that the field $\widehat {K(q)} = \widehat K(q)
_{\rm ins}$ is perfect. It follows from Lemma \ref{lemm6.1} and
\cite[Lemma~2.2]{Ch8} that in case char$(K) = 0$ (and $v$ is
discrete), one may assume without loss of generality that $v(q)
\in qv(K)$. Denote by $\Sigma $ the set of those fields $U \in
I(K(q)/K)$, for which $v(U) = v(K)$, $\widehat U \neq \widehat K$
and $\widehat U/\widehat K$ is a purely inseparable extension. In
view of Lemma \ref{lemm6.1}, our extra hypothesis ensures that
$\Sigma \neq \emptyset $. Note also that $\Sigma $ is a partially
ordered set with respect to set-theoretic inclusion, so it follows
from Zorn's lemma that it contains a maximal element, say $U
^{\prime }$. Using again Lemma \ref{lemm6.1}, one obtains that
$\widehat U ^{\prime }$ is a perfect field. Since $(K(q), v
_{K(q)})/(U ^{\prime }, v _{U'})$ is a valued extension, whence,
$\widehat {K(q)}/\widehat {U'}$ is an algebraic extension, this
proves that $\widehat {K(q)}$ is perfect as well.
\par
\smallskip
It remains to be seen that abrd$_{q}(K(q)) = 0$. Suppose first
that
\par\noindent
char$(K) = q$, fix an algebraic closure $\overline K$ of $K _{\rm
sep}$, and put $\bar v = v _{\overline K}$. It is known (cf.
\cite[Ch. VII, Theorem~22]{A1}) that if $K$ is perfect, then Br$(K
^{\prime }) _{q} = \{0\}$, for every finite extension $K ^{\prime
}/K$. We assume further that $K \neq K ^{q}$ and $K _{\rm ins}$ is
the perfect closure of $K$ in $\overline K$. It is easily verified
that $K _{\rm ins}$ equals the union $\cup _{\nu =1} ^{\infty } K
^{q^{-\nu }}$ of the fields $K ^{q^{-\nu }} = \{\beta \in
\overline K\colon \beta ^{q^{\nu }} \in K\}$, $\nu \in
\mathbb{N}$, and $[K ^{q^{-\nu }}\colon K] \ge q ^{\nu }$, for
each index $\nu $. To prove the equality abrd$_{q}(K(q)) = 0$ it
suffices to show that Br$(L ^{\prime }) _{q} = \{0\}$, for an
arbitrary $L ^{\prime } \in {\rm Fe}(K(q))$. Clearly, Br$(L
^{\prime }) _{q}$ coincides with the union of the images of Br$(L
_{0} ^{\prime }) _{q}$ under the scalar extension maps Br$(L _{0}
^{\prime }) \to {\rm Br}(L ^{\prime })$, where $L _{0} ^{\prime }$
runs across the set of finite extensions of $K$ in $L ^{\prime }$.
Moreover, one may restrict to the set $\mathcal{L}$ of those
finite extensions $L _{0} ^{\prime }$ of $K$ in $L ^{\prime }$,
for which $L _{0}'.K(q) = L ^{\prime }$ (evidently, $\mathcal{L}
\neq \emptyset $). These observations, together with basic results
on tensor products (referred to in the proof of Lemma
\ref{lemm5.7}), indicate that the concluding assertion of Lemma
\ref{lemm6.3} can be deduced from the following statement:
\par
\medskip\noindent
(6.1) Br$(L) _{q} = {\rm Br}(L.K(q)/L)$, for an arbitrary $L \in
{\rm Fe}(K)$.
\par
\medskip\noindent We prove (6.1) by showing that, for any fixed
$L$-algebra $D \in d(L)$ of
\par\noindent
$q$-power degree, there is a finite extension $K _{1}$ of $K$ in
$K(q)$ (depending on $D$), such that $[D] \in {\rm Br}(LK
_{1}/L)$, i.e. the compositum $LK _{1}$ is a splitting field of
$D$. Our proof relies on the fact that $K _{\rm ins}$ is perfect.
This ensures that Br$(K _{\rm ins} ^{\prime }) _{q} = \{0\}$
whenever $K _{\rm ins} ^{\prime } \in I(\overline K/K _{\rm
ins})$, which implies Br$(L _{1}) _{q} = {\rm Br}(L _{1}K _{\rm
ins}/L _{1})$, for every $L _{1} \in {\rm Fe}(L)$. Thus it turns
out that $[D] \in {\rm Br}(L.J'/L)$, for some finite extension $J
^{\prime }$ of $K$ in $K _{\rm ins}$. In particular, $J'$ lies in
the set, say $\mathcal{D}$, of those finite extensions $I ^{\prime
}$ of $K$ in $K _{\rm ins}$, for which $K$ has a finite extension
$\Lambda _{I'}$ in $K(q)$, such that $[D \otimes _{L} L\Lambda
_{I'}] \in {\rm Br}(L\Lambda _{I'}I'/L\Lambda _{I'})$. Choose $J
\in \mathcal{D}$ to be of minimal degree over $K$. We prove that
$J = K$, by assuming the opposite. For this purpose, we use the
following fact:
\par
\medskip\noindent
(6.2) For each $\beta \in K _{\rm ins} ^{\ast }$ and any nonzero
element $\pi \in \mathcal{M} _{v}(K(q))$, there exists $\beta
^{\prime } \in K(q) ^{\ast }$, such that $\bar v(\beta ^{\prime }
- \beta ) > v(\pi )$.
\par
\medskip\noindent
To prove (6.2) it is clearly sufficient to consider only the
special case of $\bar v(\beta ) \ge 0$. Note also that if $\beta
\in K$, then one may put $\beta ^{\prime } = \beta (1 + \pi
^{2})$, so we assume further that $\beta \notin K$. A standard
inductive argument leads to the conclusion that, one may assume,
for our proof, that $[K(\beta )\colon K] = q ^{n}$ and the
assertion of (6.2) holds for any pair $\beta _{1} \in K _{\rm ins}
^{\ast }$, $\pi _{1} \in \mathcal{M} _{v}(K(q)) \setminus \{0\}$
satisfying $[K(\beta _{1})\colon K] < q ^{n}$. Since $[K(\beta
^{q})\colon K] = q ^{n-1}$, our extra hypothesis ensures the
existence of an element $\tilde \beta \in K(q)$ with $\bar
v(\tilde \beta - \beta ^{q}) > qv(\pi )$. Applying
Artin-Schreier's theorem to the polynomial $X ^{q} - X - \tilde
\beta \pi ^{-q^{3}}$, one proves that the polynomial $X ^{q} - \pi
^{q^{2}(q-1)}X - \tilde \beta \in K(q)[X]$ has a root $\beta
^{\prime } \in K(q)$. In view of the inequality $\bar v(\tilde
\beta ) \ge 0$, this implies consecutively that $\bar v(\beta
^{\prime }) \ge 0 $ and $\bar v(\beta ^{\prime q} - \tilde \beta )
\ge q ^{2}(q-1).v(\pi )$. As $\bar v(\tilde \beta - \beta ^{q}) >
qv(\pi )$, it is now easy to see that $\bar v(\beta ^{\prime q} -
\beta ^{q}) > qv(\pi )$, whence, $\bar v(\beta ^{\prime } - \beta
) > v(\pi )$, as claimed by (6.2).
\par
\smallskip
We continue with the proof of (6.1).  The assumption that $J \neq
K$ shows that there exists $I \in I(J/K)$ with $[I\colon K] =
[J\colon K]/q$; this means that $I \notin \mathcal{D}$. Take an
element $b \in I$ so that $J = I(\sqrt[q]{b})$ and $\bar v(b) \ge
0$, and put $\Lambda = \Lambda _{J}.I$, $\Lambda ^{\prime } =
L\Lambda $. As $\widehat K(q)$ is a perfect field (i.e. $\widehat
K(q) ^{q} = \widehat K(q)$) and $v(K(q)) = qv(K(q))$, one may
assume, for our proof, that $\Lambda _{J}$ is chosen so that $b =
b _{1} ^{q}.\tilde b$, for some $b _{1} \in O _{v}(\Lambda _{J})$
and $\tilde b \in \nabla _{0}(\Lambda _{J})$.
\par
Let now $\Delta $ be the underlying division $\Lambda ^{\prime
}$-algebra of $D \otimes _{L} \Lambda ^{\prime }$. Then follows
from the Corollary in \cite[Sect. 13.4]{P}, and the choice of $J$
that $\Delta \neq \Lambda ^{\prime }$ and $[\Delta ] \in {\rm
Br}(\Lambda ^{\prime }J/\Lambda ^{\prime })$. This implies $\Delta
\cong [a, b) _{\Lambda '}$ as $\Lambda ^{\prime }$-algebras, for
some $a \in \Lambda ^{\prime \ast }$ (see, for instance, the end
of the proof of \cite[Theorem~3.2.1]{He}). It is therefore clear
that the polynomial $h _{a}(X) = X ^{q} - X - a \in \Lambda
^{\prime }[X]$ has no root in $\Lambda ^{\prime }$, so it follows
from the Artin-Schreier theorem that $h _{a}$ is irreducible over
$\Lambda ^{\prime }$, and the field $W _{a} = \Lambda ^{\prime
}(\xi _{a})$ is a degree $q$ cyclic extension of $\Lambda ^{\prime
}$, where $\xi _{a} \in \overline K$ and $h _{a}(\xi _{a}) = 0$.
One also sees that $W _{a}$ is embeddable in $\Delta $ as a
$\Lambda ^{\prime }$-subalgebra, and $\Delta $ is isomorphic to
the cyclic $\Lambda ^{\prime }$-algebra $(W _{a}/\Lambda ^{\prime
}, \sigma , b)$, for a suitably chosen generator $\sigma $ of
$\mathcal{G}(W _{a}/\Lambda ^{\prime })$. Because of the
above-noted presentation $b = b _{1} ^{q}\tilde b$, this indicates
that $\Delta \cong (W _{a}/\Lambda ^{\prime }, \sigma , \tilde
b)$. Note further that the extension $W _{a}/\Lambda ^{\prime }$
is not inertial. Assuming the opposite, one obtains from Lemma
\ref{lemm6.2} that $\tilde b \in N(W _{a}/K)$ which means that
$[\Delta ] = 0$ (cf. \cite[Sect. 15.1, Proposition~b]{P}). Since
$\Delta \in d(\Lambda ^{\prime })$ and $\Delta \neq \Lambda
^{\prime }$, this is a contradiction, proving our assertion. In
view of Ostrowski's theorem and the equality $[W _{a}\colon\Lambda
^{\prime }] = q$, the considered assertion can be restated by
saying that $\widehat W _{a}/\widehat \Lambda ^{\prime }$ is a
purely inseparable extension of degree $q$ unless $\widehat W _{a}
= \widehat \Lambda ^{\prime }$.
\par
Next we observe, using (4.1) (b), that $\eta = (\xi _{a} + 1)\xi
_{a} ^{-1}$ is a primitive element of $W _{a}/\Lambda ^{\prime }$
and $\eta \in O _{v}(W _{a}) ^{\ast }$; also, we denote by $f
_{\eta }(X)$ the minimal polynomial of $\eta $ over
$\Lambda^{\prime }$, and by $D(f _{\eta })$ the discriminant of $f
_{\eta }$. It is easily
\par\vskip0.048truecm\noindent
verified that $f _{\eta }(X) \in O _{v}(\Lambda ^{\prime })[X]$,
$f _{\eta }(0) = (-1) ^{q}$, $D(f _{\eta }) \in \Lambda ^{\prime
\ast }$, and
\par\vskip0.048truecm\noindent
$\bar v(D(f _{\eta })) = q\bar v(f _{\eta } ^{\prime }(\eta ))
> 0$ (the inequality is strict, since $[W _{a}\colon K] = q$ and
\par\vskip0.048truecm\noindent
$W _{a}/\Lambda ^{\prime }$ is not inertial). Moreover, it follows
from Ostrowski's theorem that
\par\vskip0.04truecm\noindent
there exists $\pi _{0} \in O _{v}(K)$ of value $v(\pi _{0}) =
[K(D(f _{\eta }))\colon K]\bar v(D(f _{\eta }))$. Note also that
$b ^{q ^{n-1}} \in K ^{\ast }$ (whence, $q ^{n-1}\bar v(b) \in
v(K)$), put $\pi '= \pi _{0}b ^{q ^{n-1}}$, and let $b'$ be the
$q$-th root of $b$ lying in $K _{\rm ins}$. Applying (6.2) to $b'$
and $\pi '$ (which is allowed because of the inequalities $v(\pi
') \ge v(\pi _{0}) > 0$), one obtains that there is $\lambda \in
K(q) ^{\ast }$ with $\bar v(\lambda ^{q} - b) > qv(\pi ')$.
Consider now the fields $\Lambda _{J}(\lambda )$, $\Lambda
(\lambda )$ and $\Lambda ^{\prime }(\lambda )$ instead of $\Lambda
_{J}$, $\Lambda $, and $\Lambda ^{\prime }$, respectively.
Clearly, $\Lambda _{J}(\lambda )$ is a finite extension of $K$ in
$K(q)$, $\Lambda (\lambda ) = \Lambda _{J}(\lambda ).I$ and
$\Lambda ^{\prime }(\lambda ) = L.\Lambda (\lambda )$, so our
choice of $J$ indicates that one may assume, for the proof of
(6.1), that $\lambda \in \Lambda _{J}$.
\par
\smallskip
We can now rule out the possibility that $J \neq K$, by showing that
\par\noindent
$[a, b) _{\Lambda '} \notin d(\Lambda ^{\prime })$ (in
contradiction with the choice of $J$ which requires that $I \notin
\mathcal{D}$). Indeed, the norm $N _{\Lambda '}^{W _{a}}(\lambda
\eta )$ is equal to $\lambda ^{q}$, and it follows from the
equality $\pi '= \pi _{0}b ^{q ^{n-1}}$ that $v(\pi ') \ge v(\pi
_{0}) + \bar v(b)$. Thus it turns out that
$$\bar v(\lambda ^{q}b ^{-1} - 1) > qv(\pi ') - \bar v(b) > v(\pi
_{0}) \ge \bar v(D(f _{\eta })) = q\bar v(f _{\eta } ^{\prime
}(\eta )).$$ Therefore, applying (4.1) to the polynomial $f _{\eta
}(X) + (-1) ^{q}(\lambda ^{q}b ^{-1})$ and the element $\eta $,
one obtains that $\lambda ^{q}b ^{-1}$ and $b$ are contained in
$N(W _{a}/\Lambda ^{\prime })$, which means that $[a, b) _{\Lambda
'} \notin d(\Lambda ^{\prime })$, as claimed. Hence, $J = K$, and
by the definition of $\mathcal{D}$, there exists a finite
extension $\Lambda _{K}$ of $K$ in $K(q)$, such that
\par\vskip0.04truecm\noindent
$[D \otimes _{L} L\Lambda _{K}] \in {\rm Br}(L\Lambda
_{K}/L\Lambda _{K}) = \{0\}$. In other words, $[D] \in {\rm
Br}(L\Lambda _{K}/L)$, so
\par\vskip0.04truecm\noindent
(6.1) and the equality abrd$_{q}(K(q)) = 0$ are proved in case
char$(K) = q$.
\par
\smallskip
Our objective now is to prove Lemma \ref{lemm6.3} in the special
case where $v$ is discrete. Clearly, one may assume, for our
proof, that char$(K) = 0$. Note that there exist fields $\Psi
_{\nu } \in I(K(q)/K)$, $\nu \in \mathbb{N}$, such that $\Psi
_{\nu }/K$ is a totally ramified Galois extension with $[\Psi
_{\nu }\colon K] = q ^{\nu }$ and $\mathcal{G}(\Psi _{\nu }/K)$
abelian of period $q$, for each index $\nu $, and $\Psi _{\nu '}
\cap \Psi _{\nu ''} = K$ whenever $\nu ', \nu '' \in \mathbb{N}$
and $\nu '\neq \nu ''$. This follows from \cite[Lemma~2.2]{Ch8}
(and Galois theory, which ensures that each finite separable
extension has finitely many intermediate fields). Considering, if
necessary, $\Psi _{1}$ instead of $K$, one obtains further that it
is sufficient to prove Lemma \ref{lemm6.3} under the extra
hypothesis that $v(q) \in qv(K)$. In addition, the proof of the
$q$-divisibility of $v(K(q))$ shows that, for the proof of Lemma
\ref{lemm6.3}, one may consider only the special case where
$\widehat K$ is perfect.
\par
\smallskip
Let now $\Phi $ be a finite extension of $K$ in $K _{\rm sep}$,
and $\Omega \in d(\Phi )$ a division algebra, such that $[\Omega ]
\in {\rm Br}(\Phi ) _{q}$ and $[\Omega ] \neq 0$. We complete the
proof of Lemma \ref{lemm6.3} by showing that $[\Omega ] \in {\rm
Br}(\Psi _{\nu }\Phi /\Phi )$, for every sufficiently large $\nu
\in \mathbb{N}$. As $v$ is discrete and Henselian with $\widehat
K$ perfect, the prolongation of $v$ on $\Phi $ (denoted also by
$v$) and its residue field $\widehat \Phi $ preserve the same
properties, so it follows from the assumptions on $\Omega $ that
it is a cyclic NSR-algebra over $\Phi $, in the sense of
\cite{JW}. In other words, there exists an inertial cyclic
extension $Y$ of $\Phi $ in $K _{\rm sep}$ of degree $[Y\colon
\Phi ] = {\rm deg}(\Omega )$, as well as an element $\tilde \pi
\in \Phi ^{\ast }$ and a generator $y$ of $\mathcal{G}(Y/\Phi )$,
such that $v(\tilde \pi ) \notin qv(\Phi )$ and $\Omega $ is
isomorphic to the cyclic $\Phi $-algebra $(Y/\Phi , y, \tilde \pi
)$. It follows from Galois theory and our assumptions on the
fields $\Psi _{\nu }$, $\nu \in \mathbb{N}$, that $\Psi _{\nu }
\cap Y = K$,  for all $\nu $, with, possibly, finitely many
exceptions. Fix $\nu $ so that $\Psi _{\nu } \cap Y = K$ and
deg$(\Omega ).q ^{\mu } \le q ^{\nu }$, where $\mu $ is the
greatest integer for which $q ^{\mu } \mid [\Phi \colon K]$. Put
$\Omega _{\nu } = \Omega \otimes _{\Phi } \Psi _{\nu }\Phi $ and
denote by $v _{\nu }$ the valuation of $\Psi _{\nu }\Phi $
extending $v$. It is easily obtained from Galois theory and the
choice of $\nu $ (cf. \cite[Ch. VI, Theorem~1.12]{L}) that $\Psi
_{\nu }Y/\Psi _{\nu }\Phi $ is a cyclic extension, $[\Psi _{\nu
}Y\colon \Psi _{\nu }\Phi ] = [Y\colon \Phi ] = {\rm deg}(\Omega
)$, $y$ extends uniquely to a $\Psi _{\nu }\Phi $-automorphism $y
_{\nu }$ of $\Psi _{\nu }Y$, $y _{\nu }$ generates
$\mathcal{G}(\Psi _{\nu }Y/\Psi _{\nu }\Phi )$, and $\Omega _{\nu
}$ is isomorphic to the cyclic $\Psi _{\nu }\Phi $-algebra $(\Psi
_{\nu }Y/\Psi _{\nu }\Phi , y _{\nu }, \tilde \pi )$. Also, the
assumptions on $\Psi _{\nu }$ show that $v _{\nu }(\tilde \pi )
\in q ^{\nu -\mu }v _{\nu }(\Psi _{\nu }\Phi )$. Therefore, by the
theory of cyclic algebras (cf. \cite[Sect. 15.1]{P}), and the
divisibility deg$(\Omega ) \mid q ^{\nu -\mu }$, $\Omega _{\nu }$
is $\Psi _{\nu }\Phi $-isomorphic to $(\Psi _{\nu }Y/\Psi _{\nu
}\Phi , y _{\nu }, \lambda _{\nu })$, for some $\lambda _{\nu }
\in O _{v _{\nu }}(\Psi _{\nu }\Phi ) ^{\ast }$. Since $\widehat
K$ is perfect (that is, $\widehat K = \widehat K ^{q ^{\ell }}$,
for each $\ell \in \mathbb{N}$), a similar argument shows that
$\lambda _{\nu }$ can be chosen to be an element of $\nabla
_{0}(\Psi _{\nu }\Phi )$. Taking also into account
\par\vskip0.034truecm\noindent
that $\Psi _{\nu }Y/\Psi _{\nu }\Phi $ is inertial, one obtains
from Lemma \ref{lemm6.2} that
\par\vskip0.038truecm\noindent
$\lambda _{\nu } \in N(\Psi _{\nu }Y/\Psi _{\nu }\Phi )$. Hence, by
the cyclicity of $\Psi _{\nu }Y/\Psi _{\nu }\Phi $, $[\Omega _{\nu }]
= 0$, i.e.
\par\vskip0.034truecm\noindent
$[\Omega ] \in {\rm Br}(\Psi _{\nu }Y/\Psi _{\nu }\Phi )$. As $\Psi
_{\nu } \in I(K(q)/K)$ and $\Omega \in d(\Phi )$ represents an
arbitrary nonzero element of Br$(\Phi ) _{q}$, now it becomes clear
that
\par\noindent
Br$(\Phi ) _{q} = {\rm Br}(K(q)\Phi /\Phi )$, for each $\Phi \in
{\rm Fe}(K)$, so Lemma \ref{lemm6.3} is proved.
\par
\medskip
At the end of this section, we prove two lemmas which show that
\par\noindent
dim$(K _{\rm sol}) \le 1$ whenever $K$ is a field satisfying the
conditions of Lemma \ref{lemm3.3}.
\par
\smallskip
\begin{lemm}
\label{lemm6.4} Let $(K, v)$ be a Henselian field with {\rm
char}$(\widehat K) = q$ and {\rm dim}$(\widehat K _{\rm sol})$
$\le 1$, and in case {\rm char}$(K) = 0 < q$, let $v$ be discrete.
Then {\rm dim}$(K _{\rm sol}) \le 1$.
\end{lemm}
\par
\smallskip
\begin{proof}
For each $p \in \mathbb{P}_{\widehat K}$, fix a primitive $p$-th
root of unity $\varepsilon _{p} \in K _{\rm sep}$ and a field $T
_{p}(K) \in I(K _{\rm tr}/K)$ in accordance with Lemma
\ref{lemm4.2} (b) and (c). Note first that the compositum $T(K)$
of fields $T _{p}(K)$, $p \in \mathbb{P}_{\widehat K}$, is a
subfield of $K _{\rm sol}$. Indeed, $T _{p}(K) \in I(K(\varepsilon
_{p})(p)/K)$ whenever $p \in \mathbb{P}_{\widehat K}$, so our
assertion follows from Galois theory, the cyclicity of the
extension $K(\varepsilon _{p})/K$ (cf. \cite[Ch. VI, \S{3}]{L}),
and the closeness of the class of finite solvable groups under
taking subgroups, quotient groups and group extensions. Secondly,
Lemma \ref{lemm4.1} implies the field $K _{\rm ur} \cap K _{\rm
sol} := U$ satisfies $\widehat U = \widehat K _{\rm sol}$. Using
also the fact that dim$(\widehat K _{\rm sol}) \le 1$ and $v(T(K))
= pv(T(K))$, $p \in \mathbb{P}_{\widehat K}$, and applying
(4.2)~(a), (4.3)~(a) and \cite[Theorem~2.8]{JW}, one obtains
consecutively that if $K ^{\prime }/K _{\rm sol}$ is a finite
extension, then $v(K ^{\prime }) = pv(K ^{\prime })$, Br$(K
^{\prime }) _{p} \cong {\rm Br}(\widehat K ^{\prime }) _{p}$ and
Brd$_{p}(\widehat K ^{\prime }) = {\rm Brd}_{p}(K ^{\prime }) =
0$, for each $p \in \mathbb{P}_{\widehat K}$. This proves Lemma
\ref{lemm6.4} if $q = 0$, and when $q
> 0$, our proof is completed by applying Lemma \ref{lemm6.3}.
\end{proof}
\par
\smallskip
\begin{lemm}
\label{lemm6.5} Let $K _{m}$ be an $m$-local field with {\rm
dim}$(K _{0,{\rm sol}})$ $\le 1$, $K _{0}$ being the $m$-th
residue field of $K _{m}$. Then {\rm dim}$(K _{m,{\rm sol}}) \le
1$.
\end{lemm}
\par
\smallskip
\begin{proof}
In view of Lemma \ref{lemm6.4}, one may consider only the case of
$m \ge 2$. Denote by $K _{m-j}$ the $j$-th residue field of $K
_{m}$, for $j = 1, \dots , m$. Suppose first that char$(K _{m}) =
{\rm char}(K _{0})$. Then one obtains, using repeatedly
\par\vskip0.04truecm\noindent
Lemma \ref{lemm6.4}, that dim$(K _{m,{\rm sol}}) \le 1$, which
allows to assume, for the rest of our proof, that char$(K _{m}) =
0$ and char$(K _{0}) = q > 0$. Let $\mu $ be the maximal integer
for which char$(K _{m-\mu }) = 0$. Then $0 \le \mu < m$, char$(K
_{m-\mu -1}) = q$, and in case $\mu < m - 1$, $K _{m-\mu -1}$ is
an $(m - \mu - 1)$-local field with last residue field $K _{0}$;
also, $K _{m-\mu }$ is a $1$-local field with a residue field $K
_{m-\mu -1}$. Therefore, Lemma \ref{lemm6.4} yields dim$(K _{m-\mu
',{\rm sol}}) \le 1$, for $\mu ' = \mu , \mu + 1$. Note finally
that if $\mu > 0$, then $K _{m}$ is a $\mu $-local field with $\mu
$-th residue field $K _{m-\mu }$, and by the beginning of our
proof, we have dim$(K _{m,{\rm sol}}) \le 1$, as required.
\end{proof}
\par
\smallskip
\section{\bf Tame version of Lemma \ref{lemm3.3} for admissible
Henselian fields}
\par
\medskip
Let $(K, v)$ be a Henselian field with $\widehat K$ of arithmetic
type and char$(\widehat K) = q$, put $\mathbb{P} _{q} = \mathbb{P}
\setminus \{q\}$, and suppose that abrd$_{p}(K) < \infty $, $p \in
\mathbb{P}$, and $R$ is a central division LBD-algebra over $K$.
Our main objective in this Section is to prove a modified version
of Lemma \ref{lemm3.3}, where the fields $E _{p}$, $p \in
\mathbb{P}$, are replaced by tamely ramified extensions $V _{p}$,
$p \in \mathbb{P} _{q}$, of $K$ in $K _{\rm sep}$, chosen so as to
satisfy the following conditions, for each $p \in \mathbb{P}
_{q}$: $V _{p}$ is a $p$-splitting field of $R/K$, $[V _{p}\colon
K]$ is a $p$-power, and $V _{p} \cap K _{\rm ur} \subseteq K(p)$.
The desired modification is stated as Lemma \ref{lemm7.6}, and is
also called a tame version of Lemma \ref{lemm3.3}. Our first step
towards this goal can be formulated as follows:
\par
\medskip
\begin{lemm}
\label{lemm7.1} Let $(K, v)$ be a Henselian field and let $T/K$ be
a tamely totally ramified extension of $p$-power degree $[T\colon
K] > 1$, for some $p \in \mathbb{P}_{\widehat K}$. Then there
exists a degree $p$ extension $T _{1}$ of $K$ in $T$. Moreover, $T
_{1}/K$ is a Galois extension if and only if $K$ contains a
primitive $p$-th root of unity.
\end{lemm}
\par
\smallskip
\begin{proof}
By assumption, $v(T)/v(K)$ is an abelian $p$-group of order
$[T\colon K]$, whence, there is $\theta \in T$ with $v(\theta )
\notin v(K)$ and $pv(\theta ) \in v(K)$. Therefore, by
\cite[Lemma~A.21]{TW}, there exist elements $\theta _{1} \in T$
and $a _{1} \in K$, such that \par\noindent $v(\theta _{1}) =
v(\theta )$ and $\theta _{1} ^{p} = a _{1}$. Applying
\cite[Proposition~A.22]{TW}, one concludes that the field $T _{1}
= K(\theta _{1})$ has the properties claimed by Lemma
\ref{lemm7.1}.
\end{proof}
\par
\medskip
The fields $V _{p}(K)$, $p \in \mathbb{P}_{\widehat K}$, singled
out by the next lemma play the same role in our tame version of
Lemma \ref{lemm3.3} as the role of the maximal $p$-extensions
$K(p)$, $p \in \mathbb{P}$, in the original version of Lemma
\ref{lemm3.3}.
\par
\medskip
\begin{lemm}
\label{lemm7.2} Let $(K, v)$ be a Henselian field with {\rm
abrd}$_{p}(\widehat K(p)) = 0$, for some $p \in
\mathbb{P}_{\widehat K}$. Fix $T _{p}(K) \in I(T(K)/K)$ in
accordance with Lemma \ref{lemm4.2} {\rm (c)}, and put $K _{0}(p)
= K(p) \cap K _{\rm ur}$, $V _{p}(K) = K _{0}(p).T _{p}(K)$. Then
{\rm abrd}$_{p}(V _{p}(K)) = 0$.
\end{lemm}
\par
\medskip
\begin{proof}
It follows from Lemma \ref{lemm4.2} (b) and (c) that $v(T ^{\prime
}) = pv(T ^{\prime })$, for any $T ^{\prime } \in I(K _{\rm sep}/T
_{p}(K))$; therefore, if $D ^{\prime } \in d(T ^{\prime })$ is of
$p$-power degree $\ge p$, then it is neither totally ramified nor
NSR over $T ^{\prime }$. As $p \neq {\rm char}(\widehat K)$, this
implies in conjunction with Decomposition Lemmas~5.14 and 6.2 of
\cite{JW}, that $D ^{\prime }/T ^{\prime }$ is inertial. Thus it
turns out that Br$(\widehat T ^{\prime }) _{p}$ must be
nontrivial. Suppose now that $T ^{\prime } \in I(K _{\rm sep}/V
_{p}(K))$. Then $\widehat T ^{\prime }/\widehat K(p)$ is a
separable field extension, so the condition that
abrd$_{p}(\widehat K(p)) = 0$ requires Br$(\widehat T ^{\prime })
_{p} = \{0\}$. Now it becomes clear that Br$(T ^{\prime }) _{p} =
\{0\}$, i.e. Brd$_{p}(T ^{\prime }) = 0$. Since the field $T
^{\prime }$ is an arbitrary element of $I(K _{\rm sep}/V
_{p}(K))$, this proves Lemma \ref{lemm7.2}.
\end{proof}
\par
\medskip
The following lemma presents the main properties of finite
extensions of $K$ in $V _{p}(K)$, which are used for proving Lemma
\ref{lemm3.3}.
\par
\medskip
\begin{lemm}
\label{lemm7.3} In the setting of Lemma \ref{lemm7.2}, let $V$ be
an extension of $K$ in $V _{p}(K)$ of degree $p ^{\ell } > 1$.
Then there exist fields $\Sigma _{0}, \dots , \Sigma _{\ell } \in
I(V/K)$, such that $[\Sigma _{j}\colon K] = p ^{j}$, $j = 0, \dots
, \ell $, and $\Sigma _{j-1} \subset \Sigma _{j}$ for every index
$j > 0$.
\end{lemm}
\par
\medskip
\begin{proof}
By Lemma \ref{lemm4.1} (d), the field $K$ has an inertial
extension $V _{0}$ in $V$ with $\widehat V _{0} = \widehat V$.
Moreover, it follows from (4.2) (a) and the inequality $p \neq
{\rm char}(\widehat K)$ that $V/V _{0}$ is totally ramified.
Considering the extensions $V _{0}/K$ and $V/V _{0}$, one
concludes that it is sufficient to prove Lemma \ref{lemm7.3} in
the special case where $V _{0} = V$ or $V _{0} = K$. If $V _{0} =
V$, then our assertion follows from Lemma \ref{lemm4.1} (c),
Galois theory and the subnormality of proper subgroups of finite
$p$-groups (cf. \cite[Ch. I, \S{6}]{L}). When $V _{0} = K$, by
Lemma \ref{lemm7.2}, there is a degree $p$ extension $V _{1}$ of
$K$ in $V$. As $V/V _{1}$ is totally ramified, this allows to
complete the proof of Lemma \ref{lemm7.3} by a standard inductive
argument.
\end{proof}
\par
\medskip
\label{7.5} Lemma \ref{lemm6.3} and our next lemma characterize
the fields of arithmetic type among all fields admissible by some
of Theorems \ref{theo3.1} and \ref{theo3.2}; in particular, they
show that $m$-local fields with algebraically closed last residue
fields are fields of arithmetic type, whereas an $m$-local field
with a finite $m$-th residue field is of the noted type if and
only if $m = 1$. The two lemmas also prove that if $(K, v)$ is a
Henselian field with char$(K) = {\rm char}(\widehat K)$, then $K$
is a field of arithmetic type if it is virtually perfect,
$v(K)/pv(K)\colon p \in \mathbb{P}$, are finite groups, and the
field $\widehat K$ is of arithmetic type and contains a primitive
$p$-th root of unity, for each $p \in \mathbb{P}_{\widehat K}$.
\par
\medskip
\begin{lemm}
\label{lemm7.4} Assume that $(K, v)$ and $p$ satisfy the
conditions of Lemma \ref{lemm7.2}, $\hat \varepsilon \in \widehat
K _{\rm sep}$ is a primitive $p$-th root of unity, and $\tau (p)$
is the dimension of the group $v(K)/pv(K)$, viewed as a vector
space over the field $\mathbb{Z}/p\mathbb{Z}$. Then {\rm
abrd}$_{p}(K(p)) = 0$ unless $\hat \varepsilon \notin \widehat K$
and $\tau (p) + {\rm cd}_{p}(\mathcal{G}_{\widehat K(p)}) \ge 2$.
\end{lemm}
\par
\smallskip
\begin{proof}
It follows from (4.2) (a) and Lemma \ref{lemm4.1} that if $v(K) =
pv(K)$, then $K(p) \subseteq K _{\rm ur}$, whence, $K(p) = K
_{0}(p) = V _{p}(K)$, and by Lemma \ref{lemm7.2}, abrd$_{p}(K(p))
= 0$. This agrees with the assertion of Lemma \ref{lemm7.4} in
case $v(K) = pv(K)$, since $p \neq {\rm char}(\widehat K)$ and, by
Galois cohomology, we have abrd$_{p}(\widehat K(p)) = 0$ if and
only if cd$_{p}(\mathcal{G}_{\widehat K(p)}) \le 1$ (see
\cite[Theorem~6.1.8]{GiSz} or \cite[Ch. II, 3.1]{S1}). Therefore,
we assume in the rest of the proof that $v(K) \neq pv(K)$. Fix a
primitive $p$-th root of unity $\varepsilon \in K _{\rm sep}$, and
as in Lemma \ref{lemm7.3}, consider a finite extension $V$ of $K$
in $V _{p}(K)$. It is easily verified that $\varepsilon \in K$ if
and only if $\hat \varepsilon \in \widehat K$, and this holds if
and only if $\varepsilon \in V _{p}(K)$. Suppose that $[V\colon K]
= p ^{\ell } > 1$ and take fields $\Sigma _{j} \in I(V/K)$, $j =
0, 1, \dots , \ell $, as required by Lemma \ref{lemm7.3}.
Observing that $K(p) = K _{1}(p)$, for any $K _{1} \in I(K(p)/K)$,
and using Galois theory and the normality of maximal subgroups of
nontrivial finite $p$-groups, one obtains that $V \in I(K(p)/K)$
if and only if $\Sigma _{j}/\Sigma _{j-1}$ is a Galois extension,
for every $j > 0$. In view of Lemma \ref{lemm7.1}, this occurs if
and only if $\varepsilon \in K$ or $V \in I(K _{0}(p)/K)$. It is
now easy to see that $K(p) = V _{p}(K)$ if $\varepsilon \in K$,
and $K(p) = K _{0}(p)$, otherwise. Hence, by Lemma \ref{lemm7.2},
abrd$_{p}(K(p)) = 0$ in case $\varepsilon \in K$, as claimed by
Lemma \ref{lemm7.4}.
\par
Assume finally that $v(K) \neq pv(K)$ and $\varepsilon \notin K$
(in this case, $p > 2$). It  is easy to see that if
cd$_{p}(\mathcal{G}_{\widehat K(p)}) = 1$, then there is a finite
extension $Y$ of $K _{0}(p)$ in $K _{\rm ur}$, such that $\widehat
Y(p) \neq \widehat Y$. Therefore, there exists a degree $p$ cyclic
extension $Y ^{\prime }$ of $Y$ in $Y _{\rm ur} = Y.K _{\rm ur}$,
which ensures the existence of a nicely semiramified $Y$-algebra
$\Lambda \in d(Y)$ of degree $p$; this yields
\par\vskip0.04truecm\noindent
abrd$_{p}(K _{0}(p)) \ge {\rm Brd}_{p}(Y) \ge 1$. The inequality
abrd$_{p}(K _{0}(p)) \ge 1$ also holds if $\tau (p) \ge 2$, i.e.
$v(K)/pv(K)$ is noncyclic. Indeed, then Brd$_{p}(K
_{0}(p)(\varepsilon )) \ge 1$; this follows from the fact that
$v(K _{0}(p)(\varepsilon )) = v(K)$, which implies the symbol $K
_{0}(p)(\varepsilon )$-algebra $A _{\varepsilon }(a _{1}, a _{2};
K _{0}(p)(\varepsilon ))$ (defined, e.g., in \cite{Mat}) is a
division one whenever $a _{1}$ and $a _{2}$ are elements of $K
^{\ast }$ chosen so that the cosets \par\noindent $v(a _{i}) +
pv(K)$, $i = 1, 2$, generate a subgroup of $v(K)/pv(K)$ of order
$p ^{2}$.
\par
In order to complete the proof of Lemma \ref{lemm7.4} it remains
to be seen that abrd$_{p}(K _{0}(p)) = 0$ in case
cd$_{p}(\mathcal{G}_{\widehat K(p)}) = 0$ and $\tau (p) = 1$.
Since $p \neq {\rm char}(\widehat K)$, this is the same as to
prove that cd$_{p}(\mathcal{G}_{K _{0}(p)}) \le 1$. As $K _{0}(p)
= K _{\rm ur} \cap K(p)$, we have $v(K _{0}(p)) = v(K)$ and
$\widehat {K _{0}(p)} = \widehat K(p)$, so it follows from
\cite[Lemma~1.2]{Ch4} that cd$_{p}(\mathcal{G}_{K _{0}(p)}) = {\rm
cd}_{p}(\mathcal{G}_{\widehat K(p)}) + \tau (p) = 1$, as claimed.
\end{proof}
\par
\smallskip
\begin{rema}
\label{rema7.5} Summing-up Lemmas \ref{lemm4.4}, \ref{lemm4.5} and
\ref{lemm7.4}, one obtains a complete valuation-theoretic
characterization of the fields of arithmetic type among the
maximally complete fields $(K, v)$ with abrd$_{p}(K) < \infty $,
for every $p \in \mathbb{P}$. As demonstrated in the proof of
Corollary \ref{coro4.6}, this fully describes the class
$\mathcal{C} _{0}$ of those fields of arithmetic type, which lie
in the class $\mathcal{C}$ of generalized formal power series
fields of finite absolute Brauer $p$-dimensions. Note that
$\mathcal{C}$ is considerably larger than $\mathcal{C} _{0}$. For
example, if $K _{0}$ is a finite field and $\Gamma $ is an ordered
abelian group with finite quotients $\Gamma /p\Gamma $, for all $p
\in \mathbb{P}$, then $K _{0}((\Gamma )) \in \mathcal{C} \setminus
\mathcal{C} _{0}$ in case $\Gamma /p\Gamma $ are noncyclic, for
infinitely many $p$.
\end{rema}
\par
\medskip
\label{tame} The conclusion of Lemma \ref{lemm7.3} remains valid
if $K$ is an arbitrary field, $p \in \mathbb{P}$, and $V$ is a
finite extension of $K$ in $K(p)$ of degree $p ^{\ell } > 1$; then
the extensions $\Sigma _{j}/\Sigma _{j-1}$, $j = 1, \dots , \ell
$, are Galois of degree $p$ (see \cite[Ch. I, \S{6}]{L}).
Considering the proof of Lemma \ref{lemm7.4}, one also sees that,
in the setting of Lemma \ref{lemm7.3}, $V _{p}(K) \subseteq K(p)$
if and only $\widehat K$ contains a primitive $p$-th root of unity
or $v(K) = pv(K)$. These observations and Lemma \ref{lemm7.4}
allow to view the following result as a tame version of Lemma
\ref{lemm3.3}:
\par
\medskip
\begin{lemm}
\label{lemm7.6} Assume that $K$, $q$ and $R$ satisfy the
conditions of some of Theorems \ref{theo3.1} or \ref{theo3.2}. Put
$\mathbb{P}_{q} = \mathbb{P} \setminus \{q\}$, and for each $p \in
\mathbb{P}_{q}$, denote by $k(p)$ the $p$-power of $R/K$, and by
$V _{p}(K)$ the extension of $K$ in $K _{\rm sep}$ defined in
Lemma \ref{lemm7.2}. Then there exist finite extensions $V
_{p}\colon p \in \mathbb{P}_{q}$, of $K$ in $K _{\rm sep}$ with
the following properties, for each $p \in \mathbb{P}_{q}$:
\par
{\rm (c)} $V _{p}$ is a $p$-splitting field of $R/K$, i.e. $p$ does
not divide $[V _{p}(\delta _{p})\colon V _{p}]$, for any element
$\delta _{p}$ of the underlying central division $V _{p}$-algebra
$\Delta _{p}$ of $R \otimes _{K} V _{p}$;
\par
{\rm (cc)} $V _{p} \in I(V _{p}(K)/K)$, so $[V _{p}\colon K] = p
^{\ell (p)}$, for some integer $\ell (p) \ge k(p)$, and the
maximal inertial extension $U _{p}$ of $K$ in $V _{p}$ is a
subfield of $K(p)$.
\end{lemm}
\par
\smallskip
\begin{proof}
It is clearly sufficient to show that $R \otimes _{K} V _{p}(K)
\cong M _{p^{k(p)}}(R ^{\prime })$, for an arbitrary fixed $p \in
\mathbb{P}_{\widehat K}$ and some central division $V
_{p}(K)$-algebra $R ^{\prime }$. Our proof relies on the inclusion
$V _{p}(K) \subseteq K _{\rm sol}$ and the $V _{p}(K)$-isomorphism
\par\vskip0.04truecm\noindent
$R \otimes _{K} V _{p}(K) \cong (R \otimes _{K} Y _{p}) \otimes
_{Y_{p}} V _{p}(K)$, for each $Y _{p} \in I(V _{p}(K)/K)$, which
\par\vskip0.04truecm\noindent
enables one to obtain from Lemma \ref{lemm5.1} that $R \otimes
_{K} V _{p}(K)$ is $V _{p}(K)$-isomorphic
\par\vskip0.04truecm\noindent
to $M _{s(p)}(R _{p})$, for some central division $V
_{p}(K)$-algebra $R _{p}$ and some $s(p) \in \mathbb{N}$ dividing
$p ^{k(p)}$. As a final step towards the proof of Lemma
\ref{lemm7.6}, we show that $p ^{k(p)} \mid s(p)$. Since, by Lemma
\ref{lemm7.2}, abrd$_{p}(V _{p}(K)) = 0$, it can be deduced from
\cite[Lemma~3.6]{Ch2} that for any finite extension $Y ^{\prime }$
of $V _{p}(K)$, $R _{p} \otimes _{V _{p}(K)} Y ^{\prime }$ is
isomorphic as an $Y ^{\prime }$-algebra to $M _{y'}(R ^{\prime
})$, for some $y' \in \mathbb{N}$ not divisible by $p$, and some
central division LBD-algebra $R ^{\prime }$ over $Y ^{\prime }$.
Note further that there
\par\vskip0.04truecm\noindent
is an $Y ^{\prime }$-isomorphism $R \otimes _{K} Y ^{\prime }
\cong (R \otimes _{K} Y) \otimes _{Y} Y ^{\prime }$, for any $Y
\in I(Y ^{\prime }/K)$.
\par\vskip0.04truecm\noindent
This, applied to the case of $Y = V _{p}(K)$, and combined with
the Wedderburn-Artin theorem and \cite[Sect. 9.3, Corollary~b]{P},
leads to the conclusion that $R \otimes _{K} Y ^{\prime } \cong M
_{s(p).y'}(R ^{\prime })$ as $Y ^{\prime }$-algebras. Considering
again an arbitrary $Y \in I(Y ^{\prime }/K)$, one obtains
similarly that if $R _{Y}$ is the underlying division
\par\vskip0.04truecm\noindent
$Y$-algebra of $R \otimes _{K} Y$, then there is an
$Y$-isomorphism $R \otimes _{K} Y \cong M _{y}(R _{Y})$,
\par\vskip0.04truecm\noindent
for some $y \in \mathbb{N}$ dividing $s(p).y'$. Let now $Y
^{\prime } = V _{p}(K)Y$, for some finite extension $Y$ of $K$ in
an algebraic closure of $V _{p}(K)$, such that $p ^{k(p)} \mid
[Y\colon K]$ and $Y$ embeds in $R$ as a $K$-subalgebra. Then, by
the previous observation, $p ^{k(p)} \mid s(p).y'$; as $p \nmid
y'$, this yields $p ^{k(p)} \mid s(p)$ and so completes our proof.
\end{proof}
\par
\medskip
Lemmas \ref{lemm7.2}, \ref{lemm7.3}, \ref{lemm7.6} and the results
of Sections 4 and 5 give us the possibility to deduce Lemma
\ref{lemm3.3} by the method of proving \cite[Lemma~8.3]{Ch2}. This
is done in two main steps. The first one is taken in the next
section, and our considerations there are technically facilitated
by Lemma \ref{lemm7.7}.
\par
\medskip
\section{\bf A special case of Lemma \ref{lemm3.3}}
\par
\medskip
Let $R$ be a central division LBD-algebra over a field $K$
satisfying the conditions of Theorem \ref{theo3.1} or Theorem
\ref{theo3.2}, and put $q = {\rm char}(K _{0})$ in the former
case, $q = {\rm char}(K)$ in the latter one. This section gives a
proof of Lemma \ref{lemm3.3} in the case where $q \nmid [K(r)\colon
K]$, for any $r \in R$. In order to achieve this goal we need the
following lemma:
\par
\smallskip
\begin{lemm}
\label{lemm8.1} Let $(K, v)$ be a Henselian field with {\rm
abrd}$_{\ell }(K) < \infty $, $\ell \in \mathbb{P}$, and $\widehat
K$ of arithmetic type, and let $R$ be a central division {\rm
LBD}-algebra over $K$. Fix a primitive $p$-th root of unity
$\varepsilon \in K _{\rm sep}$, for some $p \in
\mathbb{P}_{\widehat K}$, and suppose that {\rm dim}$(K _{\rm
sol}) \le 1$ and $R$ satisfies the following conditions:
\par
{\rm (i)} $p ^{2}$ and {\rm char}$(\widehat K)$ do not divide the
degree $[K(\delta )\colon K]$, for any $\delta \in R$;
\par
{\rm (ii)} There is a $K$-subalgebra $\Theta $ of $R$, which is a
totally ramified extension of $K$ of degree $[\Theta \colon K] =
p$.
\par\noindent
Then there exists a central $K$-subalgebra $\Delta $ of $R$, such
that {\rm deg}$(\Delta ) = p$ and $\Delta $ possesses a
$K$-subalgebra isomorphic to $\Theta $. Moreover, if $\varepsilon
\notin K$, then $\Delta $ contains as a $K$-subalgebra an inertial
cyclic extension of $K$ of degree $p$.
\end{lemm}
\par
\smallskip
\begin{proof}
As in the proof of Lemma \ref{lemm7.1}, one obtains from the
assumption on $\Theta /K$ that $\Theta = K(\xi )$, where $\xi $ is
a $p$-th root of an element $\theta \in K ^{\ast }$ of value
$v(\theta ) \notin pv(K)$. Suppose first that $\varepsilon \in K$.
Then $\Theta /K$ is a cyclic extension, so it follows from the
Skolem-Noether theorem that there exists $\eta \in R ^{\ast }$,
such that $\eta \xi \eta ^{-1} = \varepsilon \xi $. As a first
step towards our proof, we show that $\eta $ can be chosen so as
to satisfy the following:
\par
\medskip\noindent (8.1) The field extension $K(\eta ^{p})/K$ is
inertial.
\par
\medskip\noindent Put $\eta ^{p} = \rho $, $B = K(\rho )$, and $r =
[B\colon \mathcal{B}]$, where $\mathcal{B}$ is the maximal
inertial extension of $K$ in $B$. It is easily verified that $\xi
\rho = \rho \xi $. Since $\xi \eta \neq \eta \xi $ and
$\varepsilon \in K$, this means that $\eta \notin B$ and $[K(\eta
)\colon B] = p$. Observing that
\par
\vskip0.038truecm\noindent $[K(\eta )\colon K] = [K(\eta )\colon
B].[B\colon K]$, and by assumption, $p ^{2} \nmid [K(\eta )\colon
K]$, one also
\par
\vskip0.038truecm\noindent obtains that $p \nmid [B\colon K]$.
Therefore, $p \nmid r$, whence, the pairs $\xi , \eta $ and $\xi ,
\eta ^{r}$ generate the same $K$-subalgebra of $R$. Similarly,
condition (i) of Lemma \ref{lemm8.1} shows that char$(\widehat K)
\nmid r$, which leads to the conclusion that the set of those $b
\in B$, for which $v(b) \in rv(B)$ equals the inner group product
$\mathcal{B} ^{\ast }.\nabla _{0}(B)$. Since, by the Henselian
property of $(B, v _{B})$, $\nabla _{0}(B) \subset B ^{\ast pr}$,
this observation indicates that there exists a pair $\rho _{0} \in
\mathcal{B} ^{\ast }$, $\rho _{1} \in B ^{\ast }$, such that $\rho
^{r} = \rho _{0}\rho _{1} ^{pr}$. Putting $\eta _{1} = (\eta \rho
_{1} ^{-1}) ^{r.r'}$, for a fixed $r' \in \mathbb{N}$ satisfying
$r.r' \equiv 1 ({\rm mod} \ p)$, one obtains that $\eta _{1}\xi
\eta _{1} ^{-1} = \varepsilon \xi $ and $\eta _{1} ^{p} = \rho
_{0} ^{r'} \in \mathcal{B}$, which proves (8.1).
\par
Our objective now is to prove the existence of a $K$-subalgebra
$\Delta $ of $R$ with the properties required by Lemma
\ref{lemm8.1}. Let $\mathbb{P} ^{\prime } = \mathbb{P}_{\widehat
K} \setminus \{p\}$, and for each $p' \in \mathbb{P} ^{\prime }$,
take an extension $V _{p'}$ of $K$ in $K _{\rm tr}$ in accordance
with Lemma \ref{lemm7.6}, and put $U _{p'} = V _{p'} \cap K _{\rm
ur}$. Consider a sequence $\Pi _{n}$, $n \in \mathbb{N}$, of
pairwise distinct finite subsets of $\mathbb{P} ^{\prime }$, such
that $\cup _{n=1} ^{\infty } \Pi _{n} = \mathbb{P} ^{\prime }$ and
$\Pi _{n} \subset \Pi _{n+1}$, for each index $n$. Denote by $W
_{n}$ the compositum of the fields $V _{p _{n}}$, $p _{n} \in \Pi
_{n}$, and by $R _{n}$ the underlying division $W _{n}$-algebra of
$R \otimes _{K} W _{n}$, for any $n$. We show that $W _{n}$, $R
_{n}$ and the $W _{n}$-algebra $\Theta _{n} = \Theta \otimes _{K}
W _{n}$ satisfy conditions (i) and (ii) of Lemma \ref{lemm8.1}. It
is easily verified that $[W _{n}\colon K] = \prod _{p _{n} \in \Pi
_{n}} [V _{p _{n}}\colon K]$; in view of (4.2) (a), this ensures
that $\Theta _{n}/W _{n}$ is a totally ramified extension of
degree $p$. Using the fact that $R \otimes _{K} \Theta _{n}$ is
isomorphic to the $\Theta _{n}$-algebras $(R \otimes _{K} \Theta )
\otimes _{\Theta } \Theta _{n}$ and $(R \otimes _{K} W _{n})
\otimes _{W _{n}} \Theta _{n}$ (by \cite[Sect. 9.4,
Corollary~a]{P}), one obtains from \cite[Lemma~3.5]{Ch2} and the
uniqueness part of the Wedderburn-Artin theorem, that $\Theta
_{n}$ embeds in $R _{n}$ as a $W _{n}$-subalgebra. Note also that
$W _{n}$ and $R _{n}$ satisfy condition (i) of Lemma
\ref{lemm8.1}; since $p$ and char$(\widehat K)$ do not divide $[W
_{n}\colon K]$, this follows from Lemma \ref{lemm5.1} and
\cite[Lemma~3.5]{Ch2}.
\par
\smallskip
The next step towards our proof of the lemma can be stated as
follows:
\par
\medskip\noindent
(8.2) When $n$ is sufficiently large, $R _{n}$ has a $W
_{n}$-subalgebra $\Delta _{n} \in d(W _{n})$, such that
deg$(\Delta _{n}) = p$ and $\Theta _{n}$ embeds in $\Delta _{n}$
as a $W _{n}$-subalgebra.
\par
\medskip\noindent
Our proof of (8.2) relies on Lemma \ref{lemm5.1} and the choice of
the fields $W _{\nu }$, $\nu \in \mathbb{N}$, which indicate that
for any $\nu $ and each $\delta _{\nu } \in R _{\nu }$, $[W _{\nu
}(\delta _{\nu })\colon W _{\nu }]$ is not divisible by any $p
_{\nu } \in \Pi _{\nu }$. Arguing as in the proof of Lemma
\ref{lemm5.5}, given in \cite{Ch2}, one obtains from (8.1) the
existence of a finite dimensional $W _{\nu }$-subalgebra $\Lambda
_{\nu }$ of $R _{\nu }$ satisfying the following:
\par
\medskip\noindent
(8.3) (i) The centre $B _{\nu }$ of $\Lambda _{\nu }$ is an
inertial extension of $W _{\nu }$ of degree not divisible by
char$(\widehat K)$, $p$ and any $p _{n} \in \Pi _{n}$; moreover,
by (4.2) (a) and \par\noindent Lemma \ref{lemm4.1} (d), $B _{\nu }
= \mathcal{B} _{\nu }W _{\nu }$ and $[B _{\nu }\colon W _{\nu }] =
[\mathcal{B} _{\nu }\colon \mathcal{W} _{\nu }]$, where
$\mathcal{B} _{\nu }$ and $\mathcal{W} _{\nu }$ are the maximal
inertial extensions of $K$ in $B _{\nu }$ and $W _{\nu }$,
respectively.
\par
(ii) $\Lambda _{\nu }$ has degree $p$ as an algebra in $d(B _{\nu
})$, $C _{R _{\nu }}(\Lambda _{\nu })$ is a central division $B
_{\nu }$-algebra, and $C _{R _{\nu }}(\Lambda _{\nu })/B _{\nu }$
is of $p$-power zero (see Lemmas \ref{lemm5.1} and \ref{lemm5.3}).
\par
\medskip\noindent
It is easily verified that, for each $\nu \in \mathbb{N}$, the
field $\mathcal{W} _{\nu }$ defined in (8.3) (i) equals the
compositum of the fields $U _{p _{\nu }}$, $p _{\nu } \in \Pi
_{\nu }$. In view of Lemma \ref{lemm7.6} (cc), this means that the
Galois closures in $K _{\rm sep}$ over $K$ of the fields
$\mathcal{W} _{\nu }$, $\nu \in \mathbb{N}$, are finite Galois
extensions of $K$, and for any $\nu $,
$\mathcal{G}(\mathcal{W}_{\nu }/K)$ is nilpotent of order not
divisible by any $p' \in \mathbb{P} \setminus \Pi _{\nu }$.
Consider a $W _{\nu }$-isomorphic copy $B _{\nu } ^{\prime }$ of
$B _{\nu }$ in $W _{\nu ,{\rm sep}} = K _{\rm sep}$, and put
$\mathcal{B} _{\nu } ^{\prime } = B _{\nu } ^{\prime } \cap K
_{\rm ur}$ and $\mathcal{W} _{\nu '} = W _{\nu '} \cap K _{\rm
ur}$. It follows from (8.3) (i) and Lemma \ref{lemm4.1} (d) that
$B _{\nu } ^{\prime }/W _{\nu }$ and $\mathcal{B} _{\nu } ^{\prime
}/\mathcal{W} _{\nu }$ are inertial extensions, whereas $B _{\nu }
^{\prime }/\mathcal{B} _{\nu } ^{\prime }$ and $W _{\nu
}/\mathcal{W} _{\nu }$ are totally ramified extensions with $[B
_{\nu } ^{\prime }\colon \mathcal{B} _{\nu } ^{\prime }] = [W
_{\nu }\colon \mathcal{W} _{\nu }]$. This ensures that $B _{\nu }
^{\prime } = \mathcal{B} _{\nu } ^{\prime }W _{\nu }$ and
\par\vskip0.04truecm\noindent
$[B _{\nu } ^{\prime }\colon W _{\nu }] = [\mathcal{B} _{\nu }
^{\prime }\colon \mathcal{W} _{\nu }]$. Observing also that $[B
_{\nu } ^{\prime }\colon W _{\nu }] = [B _{\nu }\colon W _{\nu
}]$, applying Lemma \ref{lemm7.7} to the fields $\mathcal{W} _{\nu
}$, $\mathcal{B} _{\nu } ^{\prime }$ and $\mathcal{W} _{\nu '} = W
_{\nu '} \cap K _{\rm ur}$, for an arbitrary integer $\nu '
> \nu $, and using Lemma \ref{lemm4.1} (d), one obtains that:
\par
\medskip\noindent
(8.4) (a) For any pair of indices $\nu , \nu '\colon \nu < \nu '$, the
field extensions
\par\noindent
$(W _{\nu }.\mathcal{B} _{\nu } ^{\prime }\mathcal{W}_{\nu '})/(W
_{\nu }\mathcal{W}_{\nu '})$, $(\mathcal{B} _{\nu } ^{\prime
}\mathcal{W}_{\nu '})/\mathcal{W}_{\nu '}$, and $B _{\nu }
^{\prime }W _{\nu '}/W _{\nu '}$ are inertial; in
\par\vskip0.04truecm\noindent
addition, their degrees are equal and divide $[B _{\nu }\colon W
_{\nu }] = [\mathcal{B} _{\nu }\colon \mathcal{W} _{\nu }]$.
\par
(b) For each $\nu $, there exists an integer $m(\nu ) > 0$, such that
\par\noindent
$[B _{\nu } ^{\prime }W _{\nu '}\colon W _{\nu '}] = [B _{\nu }
^{\prime }W _{\nu '+1}\colon W _{\nu '+1}]$ in case $\nu ' - \nu \ge
m(\nu )$.
\par
\medskip\noindent
It is clear from (8.3) (i) and the assumptions on $\Pi _{\nu }$,
$\nu \in \mathbb{N}$, that there exists $\nu _{0} \in \mathbb{N}$,
such that, for each $\nu > \nu _{0}$, $[B _{\nu }\colon W _{\nu
}]$ has no prime divisor $p _{\nu } \le p$. Using (8.3) (i) and
(8.4), one obtains that for any $\nu > \nu _{0}$ and some $\xi
(\nu ) \in \mathbb{N}$, we have $\gcd \{[B _{\nu }\colon W _{\nu
}], [B _{n}\colon W _{n}]\} = 1$ whenever $n \in \mathbb{N}$ and
$n > \nu + \xi (\nu )$.
\par
\smallskip
We show that $R _{n}$ possesses a central $W _{n}$-subalgebra
$\Delta _{n}$ with the properties claimed by (8.2). Note that, by
Lemma \ref{lemm7.1}, $\Theta /K$ is a Galois extension, since
$\varepsilon \in K$, $\Theta /K$ is totally ramified and $[\Theta
\colon K] = p$. Fix a generator $\varphi $ of $\mathcal{G}(\Theta
/K)$, and for each $\xi \in \mathbb{N}$, let $\varphi _{\xi }$,
$\psi _{\xi }$ and $B _{\xi } ^{\prime }$ be the unique $W _{\xi
}$-automorphism of $\Theta _{\xi }$ extending $\varphi $, an
embedding of $B _{\xi }$ in $K _{\rm sep}$ as a $W _{\xi
}$-algebra, and the image of $B _{\xi }$ under $\psi _{\xi }$,
respectively. Clearly, $\psi _{\xi }$ gives rise to a canonical
bijection of $s(B _{\xi })$ on $s(B _{\xi } ^{\prime })$, which in
turn induces a homomorphism $\psi _{\xi }'\colon {\rm Br}(B _{\xi
}) \to {\rm Br}(B _{\xi } ^{\prime })$. Denote by $\Sigma _{\xi }$
and $\Sigma _{\xi } ^{\prime }$ the underlying division algebras
of $R _{\xi } \otimes _{W _{\xi }} B _{\xi }$ and $R _{\xi }
\otimes _{W _{\xi }} B _{\xi } ^{\prime }$, respectively, and let
$\widetilde B _{\xi }$ be a $W _{\xi }$-isomorphic copy of $B
_{\xi }$ in the matrix $W _{\xi }$-algebra $M _{b _{\xi }}(W _{\xi
})$, where $b _{\xi } = [B _{\xi }\colon W _{\xi }]$. Observing
that $M _{b _{\xi }}(R _{\xi }) \cong M _{b _{\xi }}(W _{\xi })
\otimes _{W _{\xi }} R _{\xi }$ over $W _{\xi }$, and applying the
Skolem-Noether theorem to $B _{\xi }$ and $\widetilde B _{\xi }$,
one obtains that $R _{\xi } \otimes _{W _{\xi }} B _{\xi } \cong M
_{b _{\xi }}(C _{R _{\xi }}(B _{\xi }))$ as $B _{\xi }$-algebras.
Hence, by the Wedderburn-Artin theorem, so are $\Sigma _{\xi }$
and $C _{R _{\xi }}(B _{\xi })$. These facts allow to identify the
$B _{\xi } ^{\prime }$-algebras $R _{\xi } \otimes _{W _{\xi }} B
_{\xi } ^{\prime }$ and $M _{b _{\xi }}(\Sigma _{\xi } ^{\prime
})$ and to prove the following:
\par
\medskip\noindent
(8.5) There is a $W _{\xi }$-isomorphism $\tilde \psi _{\xi
}\colon M _{b _{\xi }}(C _{R _{\xi }}(B _{\xi })) \to (R _{\xi }
\otimes _{W _{\xi }} B _{\xi } ^{\prime })$, extending $\psi _{\xi
}$ and mapping $C _{R _{\xi }}(B _{\xi })$ on $\Sigma _{\xi }
^{\prime }$. The image $\Lambda _{\xi } ^{\prime }$ of $\Lambda
_{\xi }$ under $\tilde \psi _{\xi }$ is a central $B _{\xi }
^{\prime }$-subalgebra of $\Sigma _{\xi } ^{\prime }$ of degree
$p$, such that $[\Lambda _{\xi } ^{\prime }] = \psi _{\xi
}'([\Lambda _{\xi }])$ in Br$(B _{\xi } ^{\prime })$.
\par
\medskip\noindent
Now fix a pair $\nu , n$ so that $\nu _{0} < \nu < \xi (\nu ) < n
- \nu $. Retaining notation as in (8.5), we turn to the proof of
the following assertion:
\par
\medskip\noindent
(8.6) The tensor products $\Lambda _{\nu } ^{\prime } \otimes _{B
_{\nu }'} (B _{\nu } ^{\prime }B _{n} ^{\prime })$, $(\Lambda
_{\nu } ^{\prime } \otimes _{B _{\nu }'} B _{\nu } ^{\prime }W
_{n}) \otimes _{B _{\nu }'W _{n}} (B _{\nu } ^{\prime }B _{n}
^{\prime })$ and $\Lambda _{n} ^{\prime } \otimes _{B _{n}'} (B
_{\nu } ^{\prime }B _{n} ^{\prime })$ are isomorphic central
division $B _{\nu } ^{\prime }B _{n} ^{\prime }$-algebras.
Moreover, the field $W _{n}$, its extensions $B _{n} ^{\prime }$
and $B _{\nu } ^{\prime }W _{n}$, and the algebras $\Lambda _{n}
^{\prime } \in d(B _{n} ^{\prime })$
\par\vskip0.04truecm\noindent
and $\Lambda _{\nu } ^{\prime } \otimes _{B _{\nu }'} (B _{\nu }
^{\prime }W _{n}) \in d(B _{\nu } ^{\prime }W _{n})$ are related as
in Lemma \ref{lemm2.2}.
\par
\medskip\noindent
The statement that $\Lambda _{\nu } ^{\prime } \otimes _{B _{\nu
}'} (B _{\nu } ^{\prime }B _{n} ^{\prime }) \cong (\Lambda _{\nu }
^{\prime } \otimes _{B _{\nu }'} B _{\nu } ^{\prime }W _{n})
\otimes _{B _{\nu }'W _{n}} (B _{\nu } ^{\prime }B _{n} ^{\prime
})$ as
\par\vskip0.032truecm\noindent
$B _{\nu } ^{\prime }B _{n} ^{\prime }$-algebras is a special case
of \cite[Sect. 9.4, Corollary~a]{P}, so it suffices
\par\vskip0.032truecm\noindent
to show that $\Lambda _{n} ^{\prime } \otimes _{B _{n}'} (B _{\nu
} ^{\prime }B _{n} ^{\prime }) \cong (\Lambda _{\nu } ^{\prime }
\otimes _{B _{\nu }'} B _{\nu } ^{\prime }W _{n}) \otimes _{B
_{\nu }'W _{n}} (B _{\nu } ^{\prime }B _{n} ^{\prime })$ over $B
_{\nu } ^{\prime }B _{n} ^{\prime }$,
\par\vskip0.07truecm\noindent
and $\Lambda _{n} ^{\prime } \otimes _{B _{n}'} (B _{\nu }
^{\prime }B _{n} ^{\prime }) \in d(B _{\nu } ^{\prime }B _{n}
^{\prime })$. Denote by $\Sigma _{\nu ,n}$ and $\Sigma _{\nu ,n}
^{\prime }$ the underlying
\par\vskip0.07truecm\noindent
division $B _{\nu } ^{\prime }B _{n} ^{\prime }$-algebras of
$\Sigma _{\nu } ^{\prime } \otimes _{B _{\nu }'} (B _{\nu }
^{\prime }B _{n} ^{\prime })$ and $\Sigma _{n} ^{\prime } \otimes
_{B _{n}'} (B _{n} ^{\prime }B _{\nu } ^{\prime })$, respectively.
\par\vskip0.056truecm\noindent
Using Lemma \ref{lemm5.7}, one obtains that $\Sigma _{\nu ,n}$ and
$\Sigma _{\nu ,n} ^{\prime }$ are isomorphic to the underlying
division $B _{\nu } ^{\prime }B _{n} ^{\prime }$-algebra of $R
\otimes _{K} B _{\nu } ^{\prime }B _{n} ^{\prime }$, whence, they
are central
\par\vskip0.04truecm\noindent division LBD-algebras
over $B _{\nu } ^{\prime }B _{n} ^{\prime }$. Note also that $p
\nmid [(B _{\nu } ^{\prime }B _{n} ^{\prime })\colon K]$; since
\par\vskip0.048truecm\noindent
$[B _{\nu } ^{\prime }B _{n} ^{\prime }\colon K] = [B _{\nu }
^{\prime }B _{n} ^{\prime }\colon W _{n}].[W _{n}\colon K]$, $B
_{\nu } ^{\prime }B _{n} ^{\prime } = B _{\nu } ^{\prime }W _{n}.B
_{n} ^{\prime }$ and (by (8.4))
\par\vskip0.04truecm\noindent
$[B _{\nu } ^{\prime }W _{n}\colon W _{n}] \mid [B _{\nu }
^{\prime }\colon W _{\nu }]$), the assertion follows from (8.3)
(i) and the fact
\par\vskip0.04truecm\noindent
that $p \nmid [W _{n}\colon K]$. In view of (8.3) (ii) and
\cite[Lemma~3.5]{Ch2}, this ensures that
\par\vskip0.04truecm\noindent
$\Lambda _{\nu } ^{\prime } \otimes _{B _{\nu }'} (B _{\nu }
^{\prime }B _{n} ^{\prime })$ and $\Lambda _{n} ^{\prime } \otimes
_{B _{n}'} (B _{\nu } ^{\prime }B _{n} ^{\prime })$ lie in $d(B
_{\nu } ^{\prime }B _{n} ^{\prime })$ and embed in $\Sigma _{\nu
,n} ^{\prime }$ as
\par\vskip0.04truecm\noindent
$B _{\nu } ^{\prime }B _{n} ^{\prime }$-subalgebras. At the same
time, $p \nmid [(B _{\nu } ^{\prime }B _{n} ^{\prime })\colon W
_{n}]$, so it follows from Lemma \ref{lemm5.1} that $\Sigma _{\nu
,n} ^{\prime }/(B _{\nu } ^{\prime }B _{n} ^{\prime })$ is of
$p$-power one. Now the proof of the former assertion of (8.6) is
completed by applying Lemma \ref{lemm5.7}, whereas the latter one
is obtained as a consequence of the choice of $\nu $, $n$, and the
observations made right after (8.4). Since, by (8.4) and the same
observations, $\gcd ([B _{\nu } ^{\prime }W _{n}\colon W _{n}], [B
_{n} ^{\prime }\colon W _{n}]) = 1$, (8.6) and Lemma \ref{lemm2.2}
imply that
\par
\medskip\noindent
(8.7) There exists $\Delta _{n} \in d(W _{n})$, such that $\Delta
_{n} \otimes _{W _{n}} B _{n} ^{\prime } \cong \Lambda _{n}
^{\prime }$ and
\par\vskip0.038truecm\noindent
$\Delta _{n} \otimes _{W _{n}} (B _{\nu } ^{\prime }W _{n}) \cong
\Lambda _{\nu } ^{\prime } \otimes _{B _{\nu }'} (B _{\nu }
^{\prime }W _{n})$ (over $B _{n} ^{\prime }$ and $B _{\nu }
^{\prime }W _{n}$, respectively).
\par
\medskip\noindent
It is clear from (8.7) and the $W _{n}$-isomorphism $B _{n} \cong
B _{n} ^{\prime }$ that $\Delta _{n} \otimes _{W _{n}} B _{n}$ and
$\Lambda _{n}$ are isomorphic $B _{n}$-algebras, which proves
(8.2). Putting $\mathcal{W} _{n} = K _{\rm ur} \cap W _{n}$, one
obtains from (8.1) that $\Delta _{n,0} \otimes _{\mathcal{W} _{n}}
W _{n} \cong \Delta _{n}$ as $W _{n}$-algebras, for some $\Delta
_{n,0} \in d(\mathcal{W} _{n})$. Note also that if $\mathcal{W}
_{n} ^{\prime }$ is the Galois closure of $\mathcal{W} _{n}$ in $K
_{\rm sep}$ over $K$, then $\mathcal{W} _{n} ^{\prime } \in I(K
_{\rm ur}/K)$, $\mathcal{G}(\mathcal{W} _{n} ^{\prime }/K)$ is
nilpotent, and the degrees $[\mathcal{W}_{n} ^{\prime }\colon K]$
and $[\mathcal{W}_{n}\colon K]$ share a common set of prime
divisors; in particular, $p \nmid [\mathcal{W} _{n} ^{\prime
}\colon K]$. This follows from Galois theory and the definition of
the fields $W _{n}$, $n \in \mathbb{N}$, which show that for each
$n$, $\mathcal{G}(\mathcal{W} _{n} ^{\prime }/K)$ is a direct
product of finite $p _{n}$-groups, where $p _{n}$ runs across
$\Pi _{n}$. Now the former assertion of Lemma \ref{lemm8.1} can be
deduced from Lemma \ref{lemm5.9} if $\varepsilon \in K$.
\par
\smallskip Let now $\varepsilon \notin K$, $[K(\varepsilon )\colon
K] = m$, and $R _{\varepsilon }$ be the underlying division
algebra of the central simple $K(\varepsilon )$-algebra $R \otimes
_{K} K(\varepsilon )$. Then $K(\varepsilon )/K$ is a cyclic field
extension and $m \mid p - 1$ (cf. \cite[Ch. VI, \S{3}]{L}), which
implies $\Theta (\varepsilon )/K(\varepsilon )$ is a totally
ramified Kummer extension of degree $p$. Observing also that $R
_{\varepsilon }$ is a central LBD-algebra over $K(\varepsilon )$,
one obtains that $\Theta (\varepsilon )$ embeds in $R
_{\varepsilon }$ as a $K(\varepsilon )$-subalgebra. At the same
time, it follows from Lemma \ref{lemm5.1} that the
\par\vskip0.031truecm\noindent
$p$-power $k(p)_{\varepsilon }$ of $R _{\varepsilon
}/K(\varepsilon )$ is less than $2$, i.e. $p ^{2} \nmid
[K(\varepsilon , \delta ')\colon K(\varepsilon )]$, for any
\par\vskip0.031truecm\noindent
$\delta ' \in R _{\varepsilon }$. Hence, $k(p)_{\varepsilon } =
1$, and by the already considered special case of our lemma, $R
_{\varepsilon }$ has a $K(\varepsilon )$-subalgebra $\Delta
_{\varepsilon } \in d(K(\varepsilon ))$, such that deg$(\Delta
_{\varepsilon }) = p$ and there exists a $K(\varepsilon
)$-subalgebra of $\Delta _{\varepsilon }$ isomorphic to $\Theta
(\varepsilon )$. Let now $\varphi $ be a generator of
$\mathcal{G}(K(\varepsilon )/K)$. Then $\varphi $ extends to an
automorphism $\bar \varphi $ of $R _{\varepsilon }$ (as a
$K$-algebra), so Lemma \ref{lemm5.4} ensures that $\Delta
_{\varepsilon }$ is $K(\varepsilon )$-isomorphic to its image
under $\bar \varphi $. Together with the Skolem-Noether theorem,
this shows that $\bar \varphi $ can be chosen so that $\bar
\varphi (\Delta _{\varepsilon }) = \Delta _{\varepsilon }$. As
$\gcd (m, p) = 1$ and $\Delta _{\varepsilon } \in d(K(\varepsilon
))$, now it follows from Teichm\"{u}ller's theorem that there is a
$K(\varepsilon )$-isomorphism
\par\vskip0.031truecm\noindent
$\Delta _{\varepsilon } \cong \Delta \otimes _{K} K(\varepsilon
)$, for some $\Delta \in d(K)$ with deg$(\Delta ) = p$. Finally,
it can be deduced from \cite[Lemma~3.5]{Ch2} that $\Delta $ is
isomorphic to a $K$-subalgebra of $R$, which in turn has a
$K$-subalgebra isomorphic to $\Theta $. Hence, by Albert's
criterion (see \cite[Sect. 15.3]{P}), $\Delta $ is a cyclic
$K$-algebra. Observe finally that cyclic degree $p$ extensions of
$K$ are inertial, since $p \neq {\rm char}(\widehat K)$ and
$\varepsilon \notin K$ (apply (4.2) (a) and Lemma \ref{lemm7.1}),
so Lemma \ref{lemm8.1} is proved.
\end{proof}
\par
\smallskip
The main lemma of the present section can be stated as follows:
\par
\medskip
\begin{lemm}
\label{lemm8.2} Let $(K, v)$ be a Henselian field with {\rm
abrd}$_{p}(K) < \infty $, $p \in \mathbb{P}$, {\rm dim}$(K _{\rm
sol}) \le 1$, and $\widehat K$ of arithmetic type; also, let $R$
be a central division {\rm LBD}-algebra over $K$ with {\rm
char}$(\widehat K) \nmid [K(\delta )\colon K]$, for any $\delta
\in R$. Then, for any $p \in \mathbb{P}_{\widehat K}$, $R/K$ has a
$p$-splitting field $E _{p}$ lying in $I(K(p)/K)$.
\end{lemm}
\par
\smallskip
\begin{proof}
Fix an arbitrary $p \in \mathbb{P}_{\widehat K}$, and a primitive
$p$-th root of unity $\varepsilon = \varepsilon _{p}$ in $K _{\rm
sep}$, take $T _{p}(K)$ as in Lemma \ref{lemm4.1} (c), and put $V
_{p}(K) = K _{0}(p).T _{p}(K)$, where $K _{0}(p) = K(p) \cap K
_{\rm ur}$. For each $z \in \mathbb{P}$, denote by $k(z)$ the
$z$-power of $R/K$, and let $\ell $ be the minimal integer $\ell
(p) \ge 0$, for which there exists an extension $V _{p}$ of $K$ in
$V _{p}(K)$ satisfying conditions (c) and (cc) of Lemma
\ref{lemm7.6}. As shown in the proof of Lemma \ref{lemm7.4}, $K(p)
= V _{p}(K)$ if $\varepsilon \in K$ or $v(K) = pv(K)$, and $K(p) =
K _{0}(p)$, otherwise. In the former case, $V _{p}$ clearly has
the properties claimed by Lemma \ref{lemm8.2}, so we suppose, for
the rest of our proof, that $\varepsilon \notin K$, $v(K) \neq
pv(K)$ and $V _{p}/K$ is chosen so that $[V _{p}\colon K] = p
^{\ell }$ and the ramification index $e(V _{p}/K)$ be minimal. Let
$E _{p}$ be the maximal inertial extension of $K$ in $V _{p}$. As
$p \in \mathbb{P}_{\widehat K}$, then it follows that $\widehat E
_{p} = \widehat V _{p}$ (cf. \cite{TW}, Proposition~A.17); using
also (4.1) (a), one sees that $V _{p}/E _{p}$ is totally ramified
and $[V _{p}\colon E _{p}] = e(V _{p}/K)$. Note further that $E
_{p} \subseteq K _{0}(p)$, by Lemma \ref{lemm7.6}, so it suffices
for the proof of Lemma \ref{lemm8.2} to show that $V _{p} = E
_{p}$ (i.e. $e(V _{p}/K) = 1$). Assuming the opposite and using
Lemma \ref{lemm7.3}, with its proof, one obtains that there is a
field $\Sigma \in I(V _{p}/E _{p})$, such that $[\Sigma \colon K]
= p ^{\ell -1}$.
\par
The main step towards the proof of Lemma \ref{lemm8.2} is to show
that $p$, the underlying division $\Sigma $-algebra $R _{\Sigma }$
of $R \otimes _{K} \Sigma $, and the field extension $V
_{p}/\Sigma $ satisfy the conditions of Lemma \ref{lemm8.1}. Our
argument relies on the assumption that dim$(K _{\rm sol}) \le 1$.
In view of Lemma \ref{lemm5.1}, it guarantees that, for each $z
\in \mathbb{P} \setminus \{p\}$, $k(z)$ is the $z$-power of $R
_{\Sigma }/\Sigma $. Thus it turns out that char$(\widehat K)
\nmid [\Sigma (\rho ')\colon \Sigma ]$, for any $\rho ' \in R
_{\Sigma }$. At the same time, it follows from the
Wedderburn-Artin theorem and the choice of $V _{p}$ and $\Sigma $
that there exist isomorphisms $R \otimes _{K} \Sigma \cong M
_{\gamma }(R _{\Sigma })$ and $R \otimes _{K} V _{p} \cong M
_{\gamma '}(R _{V_{p}})$ (as algebras over $\Sigma $ and $V _{p}$,
respectively), where $\gamma '= p ^{k(p)}$, $\gamma \mid p
^{k(p)-1}$ and $R _{V_{p}}$ is the underlying division $V
_{p}$-algebra of $R \otimes _{K} V _{p}$. Note further that the
$\Sigma $-algebras $M _{\gamma }(R _{\Sigma })$ and $M _{\gamma
}(\Sigma ) \otimes _{\Sigma } R _{\Sigma }$ are isomorphic, which
enables one to deduce from the existence of a $V _{p}$-isomorphism
$R \otimes _{K} V _{p} \cong (R \otimes _{K} \Sigma ) \otimes
_{\Sigma } V _{p}$ that $M _{\gamma '}(R _{V_{p}}) \cong M
_{\gamma }(V _{p}) \otimes _{V_{p}} (R _{\Sigma } \otimes _{\Sigma
} V _{p})$ as $V _{p}$-algebras; hence, by Wedderburn-Artin's
theorem and the inequality $\gamma < \gamma '$, $R _{\Sigma }
\otimes _{\Sigma } V _{p}$ is not a division algebra. This,
combined with \cite[Lemma~3.5]{Ch2} and the equality $[V
_{p}\colon \Sigma ] = p$, proves that $R _{\Sigma } \otimes
_{\Sigma } V _{p} \cong M _{p}(R _{V_{p}} ^{\prime })$, for some
central division $V _{p}$-algebra $R _{V_{p}} ^{\prime }$ (which
means that $V _{p}$ embeds in $R _{\Sigma }$ as a $\Sigma
$-subalgebra). It is now easy to see that
$$M _{\gamma '}(R _{V_{p}}) \cong M _{\gamma
}(V _{p}) \otimes (M _{p}(V _{p}) \otimes _{V _{p}} R _{V_{p}}
^{\prime }) \cong (M _{\gamma }(V _{p}) \otimes _{V _{p}} M _{p}(V
_{p})) \otimes _{V_{p}} R _{V_{p}} ^{\prime }$$ $$\cong M _{\gamma
p}(V _{p}) \otimes _{V _{p}} R _{V_{p}} ^{\prime } \cong M
_{p\gamma }(R _{V _{p}} ^{\prime }).$$ Using  Wedderburn-Artin's
theorem, one obtains that $\gamma = \gamma '/p = p ^{k(p)-1}$
\par\vskip0.07truecm\noindent
and $R _{V_{p}} \cong R _{V_{p}} ^{\prime }$ over $V _{p}$.
Therefore, by Lemma \ref{lemm5.1}, $p ^{2} \nmid [\Sigma (\rho
')\colon \Sigma ]$, for any $\rho ' \in R _{\Sigma }$, which
completes the proof of the fact that $p$, $R _{\Sigma }$ and $V
_{p}/\Sigma $ satisfy the conditions of Lemma \ref{lemm8.1}.
Furthermore, it follows that a finite extension of $\Sigma $ is a
$p$-splitting field of $R _{\Sigma }/\Sigma $ if and only if it is
such a field for $R/K$.
\par
\medskip
We are now in a position to complete the proof of Lemma
\ref{lemm8.2} in case $\varepsilon \notin K$ and $v(K) \neq
pv(K)$. By Lemma \ref{lemm8.1}, $R _{\Sigma }$ possesses a central
$\Sigma $-subalgebra $\Delta $, such that deg$(\Delta ) = p$ and
$V _{p}$ is embeddable in $\Delta $ as a $\Sigma $-subalgebra;
hence,
\par\vskip0.034truecm\noindent
by \cite[Theorem~4.4.2]{He}, $R _{\Sigma } = \Delta \otimes
_{\Sigma } C(\Delta )$, where $C(\Delta )$ is the centralizer of
$\Delta $ in $R _{\Sigma }$. In addition, $C(\Delta )$ is a
central division $\Sigma $-algebra, and since $p ^{2} \nmid
[\Sigma (\rho ')\colon \Sigma ]$, for any $\rho ' \in R _{\Sigma
}$, it follows from Lemma \ref{lemm5.3} that $p \nmid [\Sigma
(c)\colon \Sigma ]$, for any $c \in C(\Delta )$. Note also that
$\varepsilon \notin \Sigma $, since $\gcd ([K(\varepsilon )\colon
K], [\Sigma \colon K]) = 1$ (whence, $K(\varepsilon ) \cap \Sigma
= K$). Therefore, Lemma \ref{lemm8.1} requires the existence of a
degree $p$ cyclic extension $\Sigma ^{\prime }$ of $\Sigma $ in $K
_{\rm sep}$, which embeds in $\Delta $ as a $\Sigma $-subalgebra
(and, by Lemma \ref{lemm7.1}, is inertial over $\Sigma $). This
implies $\Sigma ^{\prime }$ is a $p$-splitting field of $R
_{\Sigma }/\Sigma $ and $R/K$ (see Lemma \ref{lemm5.3} (c) and
\cite[Sect. 13.4, Lemma and Corollary]{P}), $[\Sigma ^{\prime
}\colon K] = p ^{\ell }$, $e(\Sigma ^{\prime }/K) = e(V
_{p}/K)/p$, and $\widehat \Sigma ^{\prime }/\widehat \Sigma $ is a
cyclic extension of degree $p$. Taking finally into consideration
that $\widehat \Sigma \in I(\widehat K(p)/\widehat K)$, and using
Lemma \ref{lemm4.1} and \cite[Proposition~A.17]{TW}, one obtains
consecutively that $\widehat \Sigma ^{\prime } \in I(\widehat
K(p)/\widehat K)$ and $E _{p}$ has a degree $p$ extension $E
^{\prime }$ in $\Sigma ^{\prime } \cap K _{0}(p)$. It is now clear
that $\Sigma ^{\prime } = E ^{\prime }\Sigma $ and $\Sigma
^{\prime } \in I(V _{p}(K)/K)$. The obtained properties of $\Sigma
^{\prime }$ show that it satisfies conditions (c) and (cc) of
Lemma \ref{lemm7.6}. As $e(\Sigma ^{\prime }/K) < e(V _{p}/K)$,
this contradicts our choice of $V _{p}$ and thereby yields $e(V
_{p}/K) = 1$, i.e. $V _{p} = E _{p}$, so Lemma \ref{lemm8.2} is
proved.
\end{proof}
\par
\medskip
\section{\bf Proof of Lemma \ref{lemm3.3} and the main results}
\par
\medskip
We begin this Section with a lemma which shows how to deduce Lemma
\ref{lemm3.3} in general from its validity in the case where $q >
0 = k(q)$ ($q$ is defined at the beginning of Section 8, and
$k(q)$ is the $q$-power of $R/K$).
\par
\medskip
\begin{lemm}
\label{lemm9.1} Let $(K, v)$ be a Henselian field with $\widehat
K$ of arithmetic type, {\rm char}$(\widehat K) = q$, {\rm dim}$(K
_{\rm sol}) \le 1$ and {\rm abrd}$_{p}(K) < \infty $, $p \in
\mathbb{P}$. Take a central division {\rm LBD}-algebra $R$ over
$K$, and in case $q
> 0$, assume that $K$ has an extension $E _{q}$ in $K(q)$ that is
a $q$-splitting field of $R/K$. Then, for each $p \in
\mathbb{P}_{\widehat K}$, there is a $p$-splitting field $E _{p}$
of $R/K$, lying in $I(K(p)/K)$.
\end{lemm}
\par
\smallskip
\begin{proof}
Our assertion is contained in Lemma \ref{lemm8.2} if $q = 0$, so
we assume that $q > 0$. Let $\mathcal{R} _{q}$ be the underlying
division $E _{q}$-algebra of $R \otimes _{K} E _{q}$, and for each
$p \in \mathbb{P}$, let $k(p)'$ be the $p$-power of $\mathcal{R}
_{q}/E _{q}$. Lemma \ref{lemm5.1} (c) and the assumption on $E
_{q}$ ensure that $k(q)' = 0$, and $k(p)'$ equals the $p$-power of
$R/K$ whenever $p \in \mathbb{P}_{\widehat K}$. Therefore, by
Lemma \ref{lemm8.2}, for each $p \in \mathbb{P}_{\widehat K}$,
there is an extension $E _{p} ^{\prime }$ of $E _{q}$ in $E
_{q}(p)$, which is a $p$-splitting field of $\mathcal{R} _{q}/E
_{q}$. This enables one to deduce from Lemmas \ref{lemm5.5} and
\ref{lemm5.6} that there exist $E _{q}$-algebras $\Delta _{p}
^{\prime } \in d(E _{q})$, $p \in \mathbb{P}_{\widehat K}$,
embeddable in $\mathcal{R} _{q}$, and such that deg$(\Delta _{p}
^{\prime }) = p ^{k(p)'}$, for every $p \in \mathbb{P}_{\widehat
K}$. Hence, by Lemma \ref{lemm5.9}, $R$ possesses central
$K$-subalgebras $\Delta _{p} \in d(K)$, $p \in
\mathbb{P}_{\widehat K}$, with $\Delta _{p} \otimes _{K} E _{q}
\cong \Delta _{p} ^{\prime }$ as $E _{q}$-algebras, for each index
$p$. Moreover, it becomes clear from Lemma \ref{lemm5.3} (c) that
a finite extension $E _{p}$ of $K$ is a $p$-splitting field of
$R/K$ if and only if $E _{p}$ is a splitting field of $\Delta
_{p}/K$. In view of Lemma \ref{lemm4.3}, this proves Lemma
\ref{lemm9.1}.
\end{proof}
\par
\smallskip
We are now prepared to complete the proof of Lemma \ref{lemm3.3}
in general, and thereby to deduce Theorems \ref{theo3.1} and
\ref{theo3.2}. If $(K, v)$ and $\widehat K$ satisfy the conditions
of Theorem \ref{theo3.2}, or else, $(K, v)$ is an HDV-field with
$\widehat K$ of arithmetic type and virtually perfect, and $R$ is
a central division LBD-algebra over $K$, then the conclusion of
Lemma \ref{lemm3.3} follows from Lemmas \ref{lemm6.3} and
\ref{lemm9.1}. Therefore, as in the last paragraph of Section 5,
it is obtained that $R$ has a central $K$-subalgebra $\widetilde
R$ subject to the restrictions of Conjecture \ref{conj2.3}. This
proves Theorem \ref{theo3.2} in general, as well as Theorem
\ref{theo3.1} in the case of $m = 1$.
\par
\smallskip
In the rest of the proof of Lemma \ref{lemm3.3}, we
assume that $m \ge 2$, $K = K _{m}$ is an $m$-local field whose
$m$-th residue field $K _{0}$ is virtually perfect of
characteristic $q$ and arithmetic type, $v$ is the standard
Henselian $\mathbb{Z} ^{m}$-valued valuation of $K _{m}$ with
$\widehat K _{m} = K _{0}$, and $K _{m-m'}$ is the $m'$-th residue
field of $K _{m}$, for $m' = 1, \dots , m$. Recall that $K
_{m-m'+1}$ is complete with respect to a discrete valuation $w
_{m'-1}$ with a residue field $K _{m-m'}$, for each $m'$, and $v$
equals the composite valuation $w _{m-1} \circ \dots \circ w
_{0}$. Considering $(K, v)$ and $(K, w _{0})$, and arguing as in
the proof of Theorem \ref{theo3.2}, one obtains that the
conclusions of Theorem \ref{theo3.1} and Lemma \ref{lemm3.3} hold
if either $q = 0$ or char$(K _{m-1}) = q > 0$. It remains for us
to prove Theorem \ref{theo3.1} under the hypothesis that $q
> 0$ and char$(K _{m-1}) = 0$ (so char$(K _{m}) = 0$). Denote by
$\mu $ the maximal index for which char$(K _{m-\mu }) = 0$, fix a
primitive $q$-th root of unity $\varepsilon \in K _{\rm sep}$, put
$K ^{\prime } = K(\varepsilon )$, $v'= v _{K'}$, and denote by $R
^{\prime }$ the underlying division $K ^{\prime }$-algebra of $R
\otimes _{K} K ^{\prime }$. It is clear from Lemma \ref{lemm7.4},
applied to $R ^{\prime }$ and $(K ^{\prime }, v')$, that $R
^{\prime }$ satisfies the condition of Lemma \ref{lemm9.1},
whence, for each $p \in \mathbb{P}$, there exists a finite
extension $E _{p} ^{\prime }$ of $K ^{\prime }$ in $K ^{\prime
}(p)$, which is a $p$-splitting field of $R ^{\prime }/K ^{\prime
}$. Similarly to the proof of Lemma \ref{lemm9.1}, this allows to
show that, for each $p \in \mathbb{P}$, $R ^{\prime }$ has a $K
^{\prime }$-subalgebra $\Delta _{p} ^{\prime } \in d(K ^{\prime
})$ of degree $p ^{k(p)'}$, where $k(p)'$ is the $p$-power of $R
^{\prime }/K ^{\prime }$. Using the fact that $K ^{\prime }/K$ is
a cyclic field extension with $[K ^{\prime }\colon K] \mid q - 1$,
and applying Lemma \ref{lemm5.9} to $K ^{\prime }$, $R ^{\prime }$
and $\Delta _{q} ^{\prime }$, one concludes that $R$ has a
$K$-subalgebra $\Delta _{q} \in d(K)$, such that $\Delta _{q}
\otimes _{K} K ^{\prime } \cong \Delta _{q} ^{\prime }$ as a $K
^{\prime }$-algebra. Obviously, deg$(\Delta _{q}) = q ^{k(q)'}$,
and it follows from Lemma \ref{lemm5.1} and the divisibility $[K
^{\prime }\colon K] \mid q - 1$ that $k(q)'$ equals the $q$-power
$k(q)$ of $R/K$. Observe now that abrd$_{q}(K _{m-\mu }(q)) \le 1$
and abrd$_{q}(K _{m-\mu }) < \infty $. As char$(K _{m-\mu -1}) =
q$, the former inequality is implied by Lemma \ref{lemm6.3} and
the fact that $(K _{m-\mu }, w _{\mu })$ is an HDV-field with a
residue field $K _{m-\mu -1}$. The latter one can be deduced from
\cite[Corollary~2.5]{PS}, since $(K _{m-\mu }, w _{\mu })$ is
complete and $[K _{m-\mu -1}\colon K _{m-\mu -1} ^{q}] < \infty $.
Note also that the composite valuation $y _{\mu } = w _{\mu -1}
\circ \dots \circ w _{0}$ of $K = K _{m}$ is Henselian with a
residue field $K _{m-\mu }$ and $y _{\mu }(K) \cong \mathbb{Z}
^{\mu }$. Hence, by Lemma \ref{lemm4.4}, abrd$_{q}(K _{m}) <
\infty $. Applying finally Lemma \ref{lemm4.3} to $\Delta _{q}/K$
and $y _{\mu }$, as well as Lemma \ref{lemm5.3} to $R$ and $\Delta
_{q}$, one concludes that $(K, v)$, $q$ and $R/K$ satisfy the
conditions of Lemma \ref{lemm9.1}. Therefore, for each $p \in
\mathbb{P}$, $K$ has a finite extension $E _{p}$ in $K(p)$, which
is a $p$-splitting field of $R/K$, as claimed by Lemma
\ref{lemm3.3}. As explained at the end of Section 5, this allows
to complete the proof of Theorem \ref{theo3.1}.
\par
\medskip
Theorems \ref{theo3.1}, \ref{theo3.2} and the observations on the
scope of Lemma \ref{lemm3.3}, made in Section 3, support the idea
of making further progress on Conjecture \ref{conj2.3}, by
establishing the applicability of Lemma \ref{lemm3.3} to more
general fields $K$ with abrd$_{p}(K) < \infty $, $p \in
\mathbb{P}$. Also, it would be of interest to know whether a proof
of Conjecture \ref{conj2.3}, for a field $E$ admissible by
Proposition \ref{prop2.1}, could lead to an answer to Question
\ref{ques2.4} (b), for central division LBD-algebras over $E$.
\par
\medskip \emph{Acknowledgement.}
I would like to thank my teachers A.A. Panchishkin and P.N.
Siderov, as well as N.I. Dubrovin, A.V. Mikhalev and V.I.
Yanchevskij for the useful discussions on a number of aspects of
valuation theory and associative simple rings concerning the topic
of this paper.
\par

\end{document}